\documentclass[11pt]{article}
\edef\qedrestoreat{\noexpand\catcode\lq\noexpand\@=\the\catcode\lq\@}
\catcode\lq\@=11

\ifx\protect\undefined\let\protect\relax\fi

\def\qed{\protect\@qed{$\qedsymbol$}}
\def\pushright{\protect\@pushright}

\def\QED{\protect\@qed{{\rm Q.E.D.}}}

\def\QEI{\protect\@qed{{\rm Q.E.I.}}}

\def\Proof{\protect\@Proof}\def\endProof{\protect\@endProof}%
\let\endproof\endProof
\def\Proofof#1{\protect\@Proofof{#1}}\def\endProofof{\protect\@endProofof}%
\let\proofof\Proofof
\let\proofend\endproofof        

\def\qedsymbol{\raisebox{-.2ex}{$\Box$}}

\def\TheWordProof{\sc Proof.}
\def\TheWordProofof#1{\sc Proof of #1.}

\def\ProofFont{}

\newif\ifAutoQED\AutoQEDfalse

\newif\ifNumberResults
\ifx\theorem@style\undefined
   \NumberResultstrue
   
\fi

\def\parag@pushright#1{{
    \parfillskip=0pt            
    \widowpenalty=10000         
    \displaywidowpenalty=10000  
    \finalhyphendemerits=0      
    \hbox@pushright             
    #1
    \par}}

\def\hbox@pushright{
    \unskip                     
    \nobreak                    
    \hfil                       
    \penalty50                  
    \hskip.2em                  
    \null                       
    \hfill                      
}%

\newif\if@qed\@qedfalse
\def\save@set@qed{\let\saved@ifqed\if@qed\global\@qedtrue}%
\def\restore@qed{\global\let\if@qed\saved@ifqed}

\def\@Proof{%
   \par\removelastskip\bigskip\penalty100
   \save@set@qed
   \noindent\ProofFont{\TheWordProof\enskip}%
}%
\def\@Proofof#1{%
   \par\removelastskip\bigskip\penalty100
   \save@set@qed
   \noindent\ProofFont{\TheWordProofof{#1}\enskip}%
}%
\def\@endProof{%
   \qed\restore@qed
   \penalty-100 \medskip
}
\def\@endProofof{%
   \qed\restore@qed
   \penalty-100 \medskip
}
\def\@qed#1{%
\if@qed                                 
     \global\@qedfalse
        \ifmmode\ifinner\pushright{#1}
        \else\eqno{\qedsymbol}\fi
        \else\pushright{#1}\fi%
\else\ifhmode\ifinner\else\par\fi\fi
\fi}
\def\@pushright#1{%
  {\ifvmode                       
       \null\hfill{#1}\par        
  \else\ifmmode\maths@pushright{\hbox{#1}}
       \else\ifinner\hbox@pushright{#1}
            \else\parag@pushright{#1}
  \fi  \fi  \fi
}}%
\def\maths@pushright#1{{%
  \ifinner
     \hbox@pushright{#1}%
  \else
     \eqno#1
     \def\]{$$\ignorespaces}
  \fi
}}%
\usepackage{amsmath,amsfonts,amsbsy,amsgen,amscd, amsxtra,amstext,amsopn,amssymb}
\usepackage{latexsym}
\usepackage{bbm}
\usepackage{theorem}
\usepackage{mathrsfs}
\usepackage{latexsym}
\usepackage{euscript}

{\theorembodyfont{\slshape}
\newtheorem{theorem}{Theorem}[section]
\newtheorem{proposition}[theorem]{Proposition}
\newtheorem{lemma}[theorem]{Lemma}
\newtheorem{corollary}[theorem]{Corollary}
}
{\theorembodyfont{\rmfamily}
\newtheorem{definition}[theorem]{Definition}

\newtheorem{remark}[theorem]{Remark}
\newtheorem{example}[theorem]{Example}

}

\def\N{\mathbb{N}}

\def\R{\mathbb{R}}

\def\E{\mathbb{E}}

\def\proj{\mathbb{P}}

\newcommand{\im}{\mathrm{im}\,}

\renewcommand{\a}{\alpha}

\renewcommand{\b}{\beta}
\renewcommand{\d}{\delta}

\newcommand{\e}{\varepsilon}
\newcommand{\g}{\gamma}

\newcommand{\s}{\sigma}

\renewcommand{\t}{\tau}

\renewcommand{\tilde}{\widetilde}

\def\Oh{{\cal O}}

\newcommand\dist{\mathrm{dist}}
\newcommand\vol{\mathrm{vol}}

\newcommand\Var{\mathrm{Var}}
\newcommand\mZ{\mathcal{Z}}
\newcommand\id{\mathrm{id}}
\newcommand\tr{\mathrm{tr}}
\newcommand\diag{\mathrm{diag}}

\newcommand\zwff{\mathrm{II}}
\newcommand\calE{\mathcal{E}}
\newcommand\NJ{\mathrm{NJ}}

\def\algorithm{\begin{center}
               \begin{minipage}{6in}
               \begin{tabbing}
               \marks}
\def\falgorithm{\end{tabbing}
                \end{minipage}
                \end{center}}
\def\marks{nn\= nn\= nn\= nn\= nn\= nn\= nn\= \kill}

\begin{document}
\begin{title}
{\Large {\bf Average volume, curvatures, and Euler characteristic of random real algebraic varieties}}
\end{title}
\author{Peter B\"urgisser\thanks{Institute of Mathematics,
University of Paderborn, D-33095 Paderborn, Germany.
E-mail: {\tt pbuerg@upb.de}.
Partially supported by DFG grant BU~1371 and Paderborn Institute for Scientific Computation (PaSCo).}
}
\date{\today}
\maketitle

\begin{abstract}
We determine the expected curvature polynomial of random real projective varieties
given as the zero set of independent random polynomials with Gaussian distribution,
whose distribution is invariant under the action of the orthogonal group.
In particular, the expected Euler characteristic of such random real projective varieties
is found. This considerably extends previously known results on the number of roots,
the volume, and the Euler characteristic of the real solution set of
random polynomial equations.
\end{abstract}

\medskip


\medskip



\section{Introduction}\label{se:intro}

\subsection{Real zeros of random polynomials}

The study of the real zeros of random polyomials started with a paper by
Bloch and P\'olya~\cite{blpo:32} who investigated polynomials in one variable
with independent random coefficients from the set $\{-1,0,1\}$.
The question was further investigated by Littlewood and Offord~\cite{liof:38,liof:39}
who estimated the average number of roots with respect to different probability distributions
of the coefficients.
The first asymptotically sharp estimate on the average number of real roots
was obtained  by Kac~\cite{kac:43}. He showed that the expected number of real roots
of a random degree $d$ polynomial is asymptotically $\frac{2}{\pi}\ln d$ if the
coefficients are independent and standard normal distributed.
Erd\"os and Offord~\cite{erof:56} obtained similar results for random coefficients
in $\{-1,1\}$. Maslova~\cite{masl:74a,masl:74b} proved that Kac's asymptotic result
in fact holds for large classes of distributions.
For more details we refer to the textbook by
Bharucha-Reid and Sambandham~\cite{bhrs:86}.
Edelman and Kostlan~\cite{edko:95} give a very nice account of this and related work.

The paper by Shub and Smale~\cite{Bez2} was a breakthrough for the study of real roots
of random systems of polynomial equations. The key point in that work is an assumption
on the underlying probability measure that is very natural from a geometric point of view,
namely invariance under the action of the orthogonal group.
This probability measure was first suggested by Kostlan~\cite{kost:93}.
By sharp contrast with Kac's result~\cite{kac:43}, the expected number of real roots
of a random degree~$d$ polynomial $f=\sum_{\a=0}^d f_\a X^\a$ turns out to be exactly
$\sqrt{d}$ if the coefficients~$f_\a$ are independent centered Gaussian random variables
with variance~${d\choose \a}$. To understand the invariance, the polynomials should be
interpreted as bivariate forms (having roots in $\proj^1$).
The resulting probability distribution on bivariate forms is invariant under the
action of the orthogonal group $O(2)$.

More generally, let $H_{d,n}$ denote the vector space of homogeneous real polynomials
of degree~$d$ in the variables $X_0,\ldots,X_n$.
Let $f= \sum_{\a} f_{\a} X_0^{\a_{0}}\cdots X_n^{\a_{n}} \in H_{d,n}$ be
such that the coefficients $f_{\a}$ are independent
centered Gaussian random variables with variance
${d\choose \a}:=\frac{d!}{\a_{0}!\cdots a_{n}!}$.
The induced probability distribution on the space of forms of degree~$d$
can be shown to be invariant under
the natural action of the orthogonal group $O(n+1)$.
We will say that such $f$ is a {\em Kostlan distributed} random polynomial.
Consider a system $f_1(x)=0,\ldots,f_n(x)=0$ of Kostlan distributed random polynomials.
Shub and Smale~\cite{Bez2} showed that its expected number of real roots
equals $\sqrt{d_1\cdots d_n}$, i.e.,
the square root of the product of the degrees of the polynomials.
This result was also found by Kostlan~\cite{kost:93} in the case where all polynomials
have the same degree.
In fact, Shub and Smale's result was a byproduct
of a probabilistic analysis of nonlinear condition numbers that control the cost of
a projective homotopy method to solve systems of polynomial equations,
see also~\cite{bcss:95}.
We remark that the above choice of invariant probability measure seems also natural from
the point of view of physics~\cite{bobl:92}.

The results by Kostlan~\cite{kost:93} and Shub and Smale~\cite{Bez2} on
the expected number of real roots in the setting of an
invariant probability measure have been extended
to multihomogeneous systems by McLennan~\cite{mcle:02} and, partially, to sparse systems
by Rojas~\cite{roja:96} and Malajovich and Rojas~\cite{maro:04}.
The work of Kostlan~\cite{kost:02} contains a classification of Gaussian invariant
random polynomials along with further results.
Recently, Aza\"{\i}s and Wschebor~\cite{azw:05} gave a new proof of the Shub-Smale theorem based
on the Rice formula from the theory of random fields.
Wschebor~\cite{wsch:05}, for the first time, analyzed the variance of the number of real roots.

\subsection{Underdetermined random polynomial systems}

Considerably less is known for the underdetermined case
$f_1(x)=0,\ldots,f_s(x)=0$ ($s<n$)
where the set of solutions
is a real algebraic variety of positive dimension.
As a measure of its size different choices come to mind. One possible choice is the volume,
which is finite when the solution set is interpreted in projective space.
Another generalization of cardinality to higher dimensional solutions sets
is the Euler characteristic. (This generalization is not only natural from the topological,
but also from the computational complexity point of view, as shown in~\cite{bucu:05}.)
Both of these measures have been considered already.

Kostlan~\cite{kost:93} showed that for a system 
of Kostlan distributed $f_i$, the volume of the projective solution set has
the expectation
\begin{equation}\label{eq:kost-vol}
 \E\,\vol(\mZ(f_1,\ldots,f_s)) = d^{s/2}\, \vol(\proj^{n-s})
\end{equation}
in the case $d_1=\ldots=d_s=d$.
(Unfortunately, Kostlan never published the proof.
A proof for the more general case with possibly different $d_i$
has been given in~\cite{prio:05}.)

Podkorytov~\cite{podk:99} considered {\em any} centered Gaussian random polynomial~$f$
that is invariant under the action of the orthogonal group
and determined the expected Euler characteristic of its zero set $\mZ(f)$
in projective space $\proj^n$.
Podkorytov showed
\begin{equation}\label{eq:podk}
 \mbox{$\E\,\chi(\mZ(f)) = I_{n}(\sqrt{\d})/I_{n}(1)$\quad for odd~$n$},
\end{equation}
where $\delta$ is the {\em parameter} of~$f$ (see Definition~\ref{def:invar-poly})
and $I_{n}(\sqrt{\d)}:= \int_0^{\sqrt{\d}} (1-x^2)^{\frac{n-1}{2}} dx$.
If $f$ is Kostlan distributed, then its parameter equals the degree~$d$.
If $n$ is even, $\mZ(f)$ is almost surely a compact odd-dimensional manifold
and therefore its Euler characteristic vanishes.
Apparently, Podkorytov was not aware of Shub and Smale's work~\cite{Bez2}.
His proof is based on some tricky application of Morse theory.
An application of Podkorytov's result can be found in~\cite{kims:04}.

\subsection{Main results}

The original motivation of the present work was to extend Podkorytov's result~\cite{podk:99}
for hypersurfaces to projective varieties of higher codimension.
However, a direct use of Morse theory as in~\cite{podk:99} does not seem feasible.
It turned out to be essential to study the more general notions of
curvature coefficients and curvature polynomial of a projective variety.
Our main result (Theorem~\ref{th:main}) considerably extends and unifies
the previously known results on the number of roots, volume,
and Euler characteristic.
It determines the expectation of the curvature polynomial
of a random projective variety under an invariant probability measure.
Before stating our result, we need to explain the notion of curvature coefficients.

In a seminal work, Weyl~\cite{weyl:39} derived a formula for the volume of the tube
$T(M,\a) := \{y\in S^n\mid \dist(y,M) \le \a\}$ of radius~$\a$
around an $m$-dimensional compact smooth submanifold $M$ of the sphere $S^n$.
Let $s:=n-m$ denote the codimension of~$M$.
Weyl~\cite{weyl:39} proved that, for sufficiently small $\a>0$,
$\vol(T(M,\a))$ can be written as a linear combination
\begin{equation}\label{eq:voltube}
 \vol(T(M,\a)) = \sum_{0\le e\le m,\,\mbox{\scriptsize $e$ even}} K_{s+e}(M) J_{n,s+e}(\a)
\end{equation}
of the functions
$J_{n,k}(\a):=\int_0^{\a} (\sin\rho)^{k-1} (\cos\rho)^{n-k}d\rho$.
The coefficients~$K_{s+e}(M)$ depend on the curvature of $M$ in $S^n$
and will be thus called the {\em curvature coefficients} of the submanifold~$M$.
The functions $J_{n,k}(\a)$ determine the volume of tubes around $S^{n-k}$,
namely
$\vol(T(S^{n-k},\a)) = \Oh_{n-k}\Oh_{k-1} J_{n,k}(\a)$.

Let $\Oh_{n-1}:= 2\pi^{n/2}/\Gamma(n/2)$ denote the
$(n-1)$-dimensional volume of the unit sphere $S^{n-1}$.
Following Nijenhuis~\cite{nije:74} we rescale the curvature coefficients by
\begin{equation}\label{eq:defmu}
 \mu_e(M) := \frac1{\Oh_{m-e}\Oh_{s+e-1}} K_{s+e}(M)
 \quad \mbox{ for $0\le e\le m$, $e$ even}
\end{equation}
and define the {\em curvature polynomial} of the compact smooth submanifold~$M$
of $S^n$ by setting
\begin{equation}\label{eq:defchernpoly}
 \mu(M;T) := \sum_{0\le e\le m,\,\mbox{\scriptsize $e$ even}} \mu_e(M) T^e .
\end{equation}
For example, a subsphere $S^m$ of $S^n$ satisfies $\mu(S^m;T)=1$.
One can show that the constant term of the curvature polynomial describes
the volume of $M$, namely
$\mu_0(M) = \Oh_m^{-1}\,\vol(M)$.

The generalized Gauss-Bonnet theorem of Allendoerfer and Weil~\cite{alwe:43} and Herglotz~\cite{herg:43}
implies that the Euler characteristic $\chi(M)$ of $M$ can be retrieved by evaluating the
curvature polynomial of~$M$ at $1$: we have $ \chi(M)=2\,\mu(M;1)$ if $m$ is even,
cf.\ Theorem~\ref{th:GB-S}.

Using the canonical isometric $2$-covering map $\pi\colon S^n\to\proj^n$, we define
the curvature polynomial of a compact smooth submanifold $M$
of real projective space~$\proj^n$ by
$\mu(M;T) := \mu(\pi^{-1}(M);T))$. It is then easy to see that
$\mu_0(M) = \vol(\proj^m)^{-1}\,\vol(M)$ and
$\chi(M)=\mu(M;1)$.

Suppose now that
$f=(f_1,\ldots,f_s)\in H_{d_1,n}\times\cdots\times H_{d_s,n}$
is a Gaussian random system of polynomials.
It can be deduced from Sard's lemma
that the hypersurfaces $\mZ(f_i)$ intersect
transversally almost surely, in which case the real projective zero set
$\mZ(f_1,\ldots,f_s)\subseteq\proj^n$
is a smooth projective variety of pure codimension~$s$, or empty.
To avoid this case distinction we define $\mu(\emptyset;T):=0$.

We can now state our main result.

\begin{theorem}\label{th:main}
Suppose that $f_1\in H_{d_1,n},\ldots,f_s\in H_{d_s,n}$  ($s\le n$) are independent
centered Gaussian random polynomials with a distribution that is invariant under the action of
the orthogonal group $O(n+1)$. Let $\d_\s$ denote the parameter of $f_\s$.
Then the expected curvature polynomial of the projective
variety $\mZ(f_1,\ldots,f_s)$ in~$\proj^n$ is determined by
$$
 \E\,\mu(\mZ(f_1,\ldots,f_s);T) \ \equiv\
 \prod_{\s=1}^s \frac{\d_\s^{1/2}}{(1-(1-\d_\s)T^2)^{1/2}} \bmod T^{n-s+1}.
$$
In particular,
$\E\,\mu_e(\mZ(f_1,\ldots,f_s))$ depends only on $\d_1,\ldots,\d_s$ and $e$.
When all parameters $\d_\s$ are equal to $\d$, we obtain
$$
 \E\,\mu(\mZ(f_1,\ldots,f_s);T) \ =\
 \d^{s/2}\sum_{k=0}^{\lfloor\frac{n-s}{2}\rfloor} C_k^{(s)} (1-\d)^{k} T^{2k},
$$
where $C_0^{(s)}=1$ and
$C_k^{(s)} = \frac{s(s+2)(s+4)\cdots (s+2k-2)}{k!\,2^k}$ for $k>0$.
\end{theorem}

For instance, the theorem implies
$\E\,\vol(\mZ(f_1,\ldots,f_s)) = (\d_1\cdots\d_s)^{1/2}\,\vol(\proj^{n-s})$,
which generalizes (\ref{eq:kost-vol}).
In the case where all parameters equal~$\d$ and $n-s$ is even,
the expected Euler characteristic satisfies
\begin{eqnarray*}
 \E\,\chi(\mZ(f_1,\ldots,f_s)) &=&
 \d^{s/2}\sum_{k=0}^{\frac{n-s}{2}} C_k^{(s)} (1-\d)^{k} \\
  &=& (-1)^{\frac{n-s}{2}}\ C_{\frac{n-s}{2}}^{(s)}\ \d^\frac{n}2 + \Oh (\d^{\frac{n}2 -1})
       \quad (\d\to\infty).
\end{eqnarray*}

\subsection{Methods of proof}\label{se:methods}

The proof of our main theorem is inspired by
Aza\"{\i}s and Wschebor's~\cite{azw:05} new proof of the Shub-Smale theorem based
on the Rice formula from the theory of random fields.
Starting from Weyl's tube formula~(\ref{eq:voltube}),
we derive a version of a ``Rice formula'' for curvature coefficients
(Theorem~\ref{th:rice}) and proceed by
a probabilistic analysis of that formula,
making heavily use of the orthogonal invariance.

In fact, for the proof of Theorem~\ref{th:main}, it is sufficient to consider
the case of one equation.
This follows by employing a version of Chern's~\cite{cher:66}
kinematic formula of integral geometry for real projective space~$\proj^n$.
According to Nijenhuis~\cite{nije:74}, the expected curvature polynomial of
the intersection of a compact smooth submanifold~$M$ of $\proj^n$
with a randomly moving compact submanifold~$N$ of $\proj^n$
is a truncated product of the curvature polynomials of $M$ and $N$,
cf.~Theorem~\ref{th:kinform-S}. Actually, this was shown by Chern and Nijenhuis
for submanifolds of Euclidean space, but it also holds for projective space,
cf.\ Santal\'o~\cite{sant:76}.
Details on this can be found in the monograph by Howard~\cite{howa:93}.

Despite the considerable simplification of the arguments
for the case of a hypersurface, we develop the proof via Weyl's tube formula and
the Rice formula in full generality. Besides some intrinsic interest in obtaining
a self-contained probabilistic proof, the main reason for doing so is that
this avenue will also allow to treat higher moments (variance), as recently done so
by Wschebor~\cite{wsch:05} for the case of a zero dimensional solution set.
We thus lay the ground for a planned future paper which will investigate under which
conditions the curvature polynomial (or Euler characteristic) is
well approximated by its expectation.

We also remark that Theorem~\ref{th:main} can be quickly derived from the knowledge
of the expected Euler characteristic of a random projective hypersurface $\mZ(f)$,
for an invariant Gaussian $f$,
as derived by Podkorytov~\cite{podk:99}, cf.~(\ref{eq:podk}).
This reduction---which we found first---is again based on the kinematic formula
and the generalized Gauss-Bonnet theorem. We present it in \S\ref{se:comm}.

Finally, we remark that there is some connection of our work to the the study of
the geometric properties of random fields.
Indeed, a random polyomial $f(X_0,\ldots,X_n)$ can be seen as a ``polynomial random field''
defined on $\R^{n+1}$ or on the sphere $S^{n}$.
A central topic in Adler's book~\cite{adl:81} is the study of the
Euler characteristic (and its variation called IG characteristic)
of the ``excursion sets''
$\{x\in D\mid f(x)\ge u\}$ where $f$~is a real valued Gaussian random field on $\R^{n+1}$,
$u\in\R$, and $D\subset\R^{n+1}$ a compact domain with smooth boundary.
Part of the interest comes from the insight that
the expected Euler characteristic is a useful approximation to the distribution
of the maximum of $f$, cf.~\cite{adl:81}.
Worsley~\cite{wors:95} describes some applications of the statistics
of the Euler characteristic of excursion sets in $\R^n$ to astrophysics and medicine.
Adler and Taylor~\cite{taad:03} have considerably extended and unified the previous results
on the expected Euler characteristic of excursion sets to the general framework of a centered
regular Gausssian field on a compact manifold.
We realized that it is possible to deduce the expected Euler characteristic of
a random projective hypersurface, that is,
Podkorytov's result~\cite{podk:99}, from the general result in~\cite[Theorem~4.1]{taad:03},
cf.\ Remark~\ref{re:adta}.
We remark that our proof is methodically quite different from the ones by
Podkorytov, and Adler and Taylor. Both of them rely on Morse theory, while
we analyze expected tube volumes.
As already mentioned before, the direct application of Morse theory does not seem
feasible for the probabilistic analysis of projective varieties of higher codimension.

The structure of the paper is roughly as follows:
In \S\ref{se:diffintgeo} we recall relevant facts from differential and integral geometry.
Sections~\ref{se:vectpolymat}--\ref{se:invarpoly} prepare for the proof of Theorem~\ref{th:main},
which is then given in \S\ref{se:randvar}.
Hereby, \S\ref{se:vectpolymat} develops the necessary facts about
Gaussian random vectors and symmetric matrices that are invariant under the action of the
orthogonal group. In \S\ref{se:invarpoly} we give a discussion on invariant random polynomials.
In \S\ref{se:comm} we show how to quickly derive Theorem~\ref{th:main}
from the knowledge of the expected Euler characteristic of a random projective hypersurface
for an invariant centered Gaussian random polynomial.


\subsection*{Acknowledgments}

I thank Martin Lotz and Mario Wschebor for useful discussions.
I am grateful to Mario Wschebor for pointing out to me
Robert Adler's book, and I thank Alexander Alldridge for finding the
reference to the monograph by Ralph Howard.
Finally, I thank Dima Grigoriev for quickly sending me
Podkorytov's paper.

\section{Background from differential and integral geometry}\label{se:diffintgeo}

\subsection{Weyl's tube formula and curvature polynomials in spheres}\label{se:weyltube}

For the following material from differential geometry we refer e.g.\
to \cite{kono2:69} or \cite{spiv3:79,spiv4:79}.
Let $M$ be a compact smooth $m$-dimensional submanifold of $S^n$,
interpreted as a Riemannian submanifold.
We denote by $T_xM$ the tangent space of $M$ at a point $x\in M$ and write $(T_xM)^\perp$
for its orthogonal complement in $T_xS^n$.
The curvature of $M$ at~$x$ is described by the {\em second fundamental form}
of $M$ at $x$, which is a trilinear map
$\zwff_M(x)\colon T_xM\times T_xM\times (T_xM)^\perp\to\R$ that is symmetric in the first
two components.
In terms of local coordinates, $\zwff_M$ can be described as follows: let
$u=(u_1,\ldots,u_{m})\mapsto \varphi(u_1,\ldots,u_{m})\in M\subset\R^{n+1}$
be a local parametrization of $M$. Then, for any unit normal vector $\nu\in (T_xM)^\perp$,
\begin{equation}\label{eq:def-2ndff}
 \zwff_M(x)\big(\partial_\a,\partial_\b,\nu\big)
  = \big\langle\partial^2_{\a,\b}\varphi,\nu\big\rangle,
\end{equation}
using the shorthand notation
$\partial_\a :=\partial_{u_\a}$ and $\partial^2_{\a,\b} :=\partial^2_{u_\a,u_\b}$.
By the {\em Weingarten map} of $M$ at $x$ in direction~$\nu$
we understand the self adjoint linear map
$L_M(x,\nu)\colon T_xM \to T_xM$
characterized by
\begin{equation}\label{eq:defoL}
\langle L_M(x,\nu)(V),W\rangle = \zwff_M(x)(V,W,\nu)\quad
\mbox{for $V,W\in T_xM$.}
\end{equation}

In a seminal work, Weyl~\cite{weyl:39} determined the volume of the {\em tube}
$$
T(M,\a) := \{y\in S^n\mid \dist(y,M) \le \a\}
$$
around~$M$ for a sufficiently small radius $\a>0$.
He proved that
\begin{equation}\label{eq:weyl}
 \vol(T(M,\a)) = \int_{x\in M} \int_{\rho=0}^a \int_{\nu\in S_x}
   \frac{\rho^{s-1}  \det(\id - \rho L_M(x,\nu))}{(1+\rho^2)^{(n+1)/2}}\, dS_x(\nu)\, d\rho\, dM(x) ,
\end{equation}
where $a=\tan\a$, $S_x$ denotes the unit sphere in $(T_xM)^\perp$,
and $s=n-m$ is the codimension of $M$ in $S^n$
(see also~\cite{gray:90}).
Moreoever, Weyl showed that the tube volume $\vol(T(M,\a))$
can be written as a linear combination
of the linearly independent functions $J_{n,s+e}(\a)$
with real coefficients~$K_{s+e}(M)$,
for $0\le e\le m$, $e$ even,
cf. Equation (\ref{eq:voltube}).
The lowest order coefficient $K_{s}(M)$ equals $\Oh_{s-1}\vol(M)$,
which is intuitively plausible.

We will call $K_{s+e}(M)$ the
{\em curvature coefficients of the submanifold~$M$ of~$S^n$}.
In order to justify the naming of these coefficients, we remark that,
after some rescaling, $K_{s+e}(M)$ is an isometric invariant of the
Riemannian submanifold $M$ of~$S^n$. More precisely, for $e>0$,
$$
 \frac{s (s+2)\cdots (s+e-2)}{\Oh_{s-1}} K_{s+e} = \int_M k_e dM
$$
with some function $k_e\colon M\to\R$ whose value at $x\in M$
depends only on the difference of the Riemann tensor of $M$
and the Riemann tensor of $S^n$ restricted to $M$, at~$x$,
cf.~\cite{weyl:39}.
(We will not need this observation in the following.)
It is important to realize that these curvature coefficients
are not ``absolute'' invariants of the Riemannian manifold~$M$.
One can show that for a subsphere $M=S^m$ of  $S^n$
we have $K_{s+e}(S^m)=0$ for $e\ne 0$.

For our purposes, it will be more useful to rescale the curvature coefficients
as done by Nijenhuis~\cite{nije:74}:
we define the {\em curvature polynomial}
$ \mu(M;T)$ 
as in the introduction (Equations~(\ref{eq:defmu}) and (\ref{eq:defchernpoly})).
It follows from the above that the constant term of the curvature polynomial describes
the volume of~$M$: we have $\mu_0(M) = \Oh_m^{-1}\,\vol(M)$.
We note that $\mu(S^m;T)=1$ for a subsphere $S^m$ of $S^n$.

The curvature polynomial also encodes the Euler characteristic of $M$ in a simple way.
The following statement can be deduced from the generalized Gauss-Bonnet theorem
of Allendoerfer and Weil~\cite{alwe:43} and
Herglotz~\cite{herg:43}. We provide the proof in the appendix.

\begin{theorem}\label{th:GB-S}
Let $M$ be a compact smooth submanifold of $S^n$ of even dimension~$m$.
Then the Euler characteristic of $M$ can be expressed as
$\chi(M)=2\,\mu(M;1)$.
\end{theorem}

\subsection{Principal kinematic formula of integral geometry for spheres}\label{se:kin-form}

One of the main goals of integral geometry is to compute
integrals of the form $\int_G I(M\cap gN) dg$,
where $M$ and $N$ are compact smooth submanifolds of a homogenous space
with respect to the action of a Lie group~$G$, $I$ is some integral invariant,
and the integration is with respect to the invariant measure on $G$.
Kinematic formulas provide answers to this question in the form
$\int_G I(M\cap gN) dg  = \sum_k c_k I_k(M) J_k(N)$ with integral invariants $I_k,J_k$ related to $I$.
For a comprehensive treatment of this subject we refer to Santal\'o's book~\cite{sant:76}.
A unified treatment of  kinematic formulas in homogeneous spaces has been given
by~Howard~\cite{howa:93}.

Chern~\cite{cher:66} and Federer~\cite{fede:59} proved a general
kinematic formula for submanifolds of Euclidean space with respect to the group of motions.
Nijenhuis~\cite{nije:74} pointed out a particular elegant formulation of this kinematic formula.
He observed that, after some rescaling of integral invariants, $\sum_k c_k I_k(M) J_k(N)$ can be
interpreted as a reduced polynomial multiplication. This leads to a great deal of simplification
in our calculations, as the formulas for the coefficients $c_k$ turn out to be quite complicated.

For submanifolds of the sphere and the orthogonal group, the kinematic formula takes exactly the
same form as for Euclidean space.
An indication of this at first glance astonishing fact can be found, somewhat hidden,
in Santal\'o~\cite[IV.18.3. p.~320]{sant:76} for the special case of the intersection of domains.
Howard~\cite{howa:93} clarified this phenomenon by establishing a general transfer theorem
according to which the Chern-Federer kinematic formulas hold in all simply connected
homogeneous spaces of constant sectional curvature and not just in Euclidean space.

The kinematic formula for spheres allows the following beautiful formulation.

\begin{theorem}\label{th:kinform-S}
Let $M$ and $N$ be compact smooth submanifolds of $S^n$
of the dimensions $m$ and $p$, respectively, such that $m+p\ge n$.
Then we have
$$
 \int \mu(M\cap gN;T) dg \equiv
  \mu(M;T) \mu(N;T) \bmod T^{m+p-n+1},
$$
where the integration is with respect to the Haar measure on the
orthogonal group $O(n+1)$ scaled such that the volume of $O(n+1)$ equals $1$.
In particular,  we have
$$
\int \mu(M\cap gS^p;T) dg
 \equiv \mu(M;T) \bmod T^{m+p-n+1}.
$$
\end{theorem}

We note that Poincar\'e's formula
$$
\int  \Oh_{m+p-n}^{-1}\vol(M\cap gN)  dg = \Oh_m^{-1} \vol(M)\  \Oh_p^{-1}\vol(N),
$$
is a special case of Theorem~\ref{th:kinform-S}, obtained
by comparing constant coefficients.

Weyl's tube formula and the kinematic formula
immediately extend from $S^n$ to the real projective space $\proj^n$,
using the canonical isometric $2$-covering map
$\pi\colon S^n\to\proj^n$.
We define the curvature polynomial of a compact submanifold~$M$
of real projective space~$\proj^n$ by
\begin{equation}\label{eq:defmuP}
 \mu(M;T) := \mu(\pi^{-1}(M);T) .
\end{equation}
This definition gives the appropriate scaling, as
$$
\mu_0(M) = \mu_0(\pi^{-1}(M)) = \frac1{\Oh_m} \vol(\pi^{-1}(M))
 =  \frac2{\Oh_m} \vol(M) = \frac1{\vol(\proj^m)} \vol(M).
$$
Moreover, it is easy to see that the  kinematic formula of Theorem~\ref{th:kinform-S}
also holds for~$\proj^n$.
Finally, as $\pi\colon S^n\to\proj^n$ is a two sheeted covering map,
Theorem~\ref{th:GB-S} implies that
$\chi(M) = \frac12 \chi(\pi^{-1}(M)) = \mu(\pi^{-1}(M);1) = \mu(M;1)$.

\begin{remark}\label{re:alt-pf-kf}
It is possible to derive Theorem~\ref{th:kinform-S} from the kinematic formulas given
in Santal\'o~\cite[(15.72), p.~269 and IV.18.3. p.~320]{sant:76}
for the intersection of domains in Euclidean space and spheres, respectively.
Santal\'o~\cite[p.~222, p.~302]{sant:76}
assigns to a compact domain $Q$ in $\R^n$ or $S^n$, bounded by a smooth hypersurface
$\partial Q$, the following integral of mean curvature
\begin{equation*}
 M_i^{San}(\partial Q) := {n-1\choose i}^{-1}
  \int_{\partial Q}\sigma_i(\kappa_1,\ldots,\kappa_{n-1})\,d(\partial Q) .
\end{equation*}
Hereby, $\sigma_i(\kappa_1,\ldots,\kappa_{n-1})$
stands for the $i$th elementary symmetric function
in the principal curvatures $\kappa_j$ of the hypersurface $\partial Q$.
The curvature coefficients $K_{i+1}(M)$ of a a compact
submanifold $M$ of $S^n$ of codimension~$s$ 
can be related to Santal\'o's integral of mean curvature
of the tubes around $M$ as follows 
\begin{equation}\label{eq:meancurv}
 K_{i+1}(M) =  {n-1\choose i} \lim_{\a\to 0} M_{i}^{San}(\partial T(M,\a)).
\end{equation}
Herebye, $K_j(M)$ is defined by (\ref{eq:voltube})  if $j\ge s$ and $j-s$ is even. 
Otherwise, we set $K_j(M)=0$.
The proof is as in \cite[V\S4]{suwi:72}, where a similar result is shown for Euclidean space.
\end{remark}

\section{Invariant Gaussian vectors and matrices}\label{se:vectpolymat}

Here we develop facts about invariant random matrices that will be used in
\S\ref{se:randvar} for the proof of the main result.

A random vector is called centered iff its expectation is zero.
We call two random vectors {\em equivalent} if they have the same distribution.
In particular, equivalent random vectors have the same expectation and covariance matrix.
We shall write $X\sim N(0,\sigma^2)$ to indicate that $X$ is a real valued Gaussian variable
with mean zero and variance $\sigma^2$.

\subsection{Invariant random vectors}\label{se:invar-vect}

A random vector $V\in\R^n$ is called {\em $O(n)$-invariant}
iff $gV$ is equivalent to $V$ for all $g\in O(n)$.
For simplicity, we assume that $V$ has a density.

\begin{lemma}\label{le:univect}
Let $V\in\R^n$ be an $O(n)$-invariant random vector. Then
$V/\|V\|$ is uniformly distributed in $S^{n-1}$ and independent of $\|V\|$.
\end{lemma}

\begin{proof}
By invariance, the density of $V$ depends only on $\|V\|$.
\end{proof}

Invariant Gaussian random vectors are easy to characterize
by the following well known fact.

\begin{lemma}\label{le:invar-vect}
A Gaussian random vector $V=(V_1,\ldots,V_n)$ is $O(n)$-invariant iff
it is centered and
$V_1,\ldots,V_n$ are independent and have the same variance.
\end{lemma}

\begin{proof}
Suppose $V$ is $O(n)$-invariant. Then $-V$ is equivalent to $V$
(take $g=-I_n$), hence
$\E\,V=0$. By assumption, $gV V^T g^T$ is equivalent to $V V^T$, hence
the covariance matrix $A:=\E(V V^T)$ satisfies $gAg^T=A$ for all $g\in O(n)$.
A straightforward calculation shows that $A$ must be a multiple of the
unit matrix $I_n$ (take for $g$ rotations in two dimensional
coordinate subspaces).
The converse follows from the fact that Gaussian random variables are
characterized by their expectation and covariance matrix.
\end{proof}

\subsection{Invariant random symmetric matrices}\label{se:invar-mat}

Let $\Sigma_n$ denote the space of real symmetric $n$ by $n$ matrices.
The Frobenius norm of $W\in\Sigma_n$ is defined as
$\|W\|_F:=(\sum_{i,j} W_{ij}^2)^{1/2}$.
We will assume $n>1$.

\begin{definition}\label{def:parW}
The {\em parameter $\delta(W)$} of a random matrix $W\in\Sigma_n$ is defined as
$$
  \delta(W):= \frac1{n(n-1)} \big(\E\,(\tr W)^2 - \E\,\|W\|_F^2 \big).
$$
A random matrix $W\in\Sigma_n$ is called {\em $O(n)$-invariant} iff
$gWg^T$ is equivalent to $W$ for all $g\in O(n)$.
\end{definition}

The following proposition classifies all invariant Gaussian $W\in\Sigma_n$,
compare~\cite[Prop.~4.1]{podk:99}.

\begin{proposition}\label{le:classW}
Suppose $T\in\Sigma_n$ with independent $T_{ij}$ and
$T_{ij}\sim N(0,1)$ for $i\ne j$ and $T_{ii}\sim N(0,2)$, $n>1$.
Moreover, let $r,s\in\R$ and $Z\sim N(0,1)$ be independent of~$T$.
Then $W=rZ I_n + sT$ is $O(n)$-invariant and has the parameter
$\delta(W)=r^2-s^2$.
Moreover, any $O(n)$-invariant Gaussian $W\in\Sigma_n$
is equivalent to one of this form, in particular $\E\,W=0$.
\end{proposition}

\begin{proof}
In order to see that $T$ is $O(n)$-invariant it is sufficient to
check that $\tilde{T}=gTg^T$ is equivalent to~$T$ for a set of generators~$g$.
Since $\tilde{T}$ and $T$ are both centered Gaussian it is sufficient to check that
they have the same covariance matrix.
The group $O(n)$ is generated by $\diag(-1,1,\ldots,1)$ and the rotations
in two dimensional coordinate spaces. Invariance under the latter is verified
by a straightforward calculation of covariances.
Moreover, a calculation shows that
$$
 \E(\tr W)^2 = n^2 r^2 + 2n s^2,\quad
 \E\,\|W\|_F^2 = n r^2 + n(n+1) s^2 .
$$
This implies $\delta(W)= r^2-s^2$.

Suppose now that $W$ is $O(n)$-invariant.
Let $g$ be the product of the permutation matrix $P_{\pi^{-1}}$ and $\diag(\e_1,\ldots,\e_n)$
where $\e_i=\pm 1$. Then $gWg^T=(\e_i\e_jW_{\pi(i),\pi(j)})$ is equivalent to $W$.
Taking $\pi=\id,\e_1=-1,\e_2=\ldots\e_n=1$ we conclude that
$(W_{11},W_{12})$ is equivalent to $(W_{11},-W_{12})$, hence
$\E\, W_{12} =0$ and
$\E\,(W_{11} W_{12}) = 0$.
By choosing appropriate $\pi$ and $\e_i$ one can similarly show that
$\E\,(W_{ij}W_{k\ell})=0$ except in the cases where $i=j$ and $k=\ell$
(details are left to the reader).
Similarly, one shows that $W$ is centered.
Moreover, by conjugating with permutation matrices, we see that
$$
\E\, W_{ii}^2=\E\, W_{11}^2,\ \E\, W_{ij}^2 = \E\, W_{12}^2 =: s^2,\
\E\,(W_{ii} W_{jj}) = \E\,(W_{11}W_{22}) =:r^2 \quad (i\ne j).
$$
If we can prove that
\begin{equation}\label{eq:covare}
 \E W_{11}^2 = \E\,(W_{11} W_{22}) + 2\, \E\, W_{12}^2 ,
\end{equation}
then $W$ has the same covariance matrix as $rZ I_n + sT$ and we are done.

For showing (\ref{eq:covare}) suppose without loss of generality $n=2$.
Put $\tilde{W}=gWg^T$ where
$g=\begin{pmatrix} \cos\varphi & -\sin\varphi\\ \sin\varphi & \cos\varphi \end{pmatrix} $.
Then
$\tilde{W}_{11} = W_{11}\cos^2\varphi + W_{22}\sin^2\varphi - 2W_{12}\cos\varphi\sin\varphi$.
Hence
$$
  \E \tilde{W}_{11}^2 = \E\, W_{11}^2 (\cos^4\varphi + \sin^4\varphi)
  + 2\,\Big(2\,\E\,W_{12}^2 + \E\,(W_{11} W_{22})\Big) \cos^2\varphi\sin^2\varphi,
$$
Using that $\cos^4\varphi + \sin^4\varphi = 1 -2\,\cos^2\varphi\sin^2\varphi$,
the assertion~(\ref{eq:covare}) follows.
\end{proof}

\begin{lemma}\label{le:suppositionW}
Suppose $W,W_1,\ldots,W_s\in \Sigma_n$ are random symmetric matrices
and $\lambda_1,\ldots,\lambda_s\in\R$.
Further, let $u$ be a real random variable. Then:
\begin{description}
\item[(i)] $\delta(W + u I_n) = \delta(W) + \E\,u^2 + \frac2{n}\E\,(u\,\tr W)$.

\item[(ii)] $\delta( \lambda_1 W_1 +\cdots+\lambda_sW_s)
 = \lambda_1^2\delta(W_1) + \cdots + \lambda_s^2\delta(W_s)$
if $W_1,\ldots,W_s$ are independent.
\end{description}
\end{lemma}

\begin{proof}
(i) Put $\tilde{W} = W + uI_n$. Then
$(\tr\tilde{W})^2 = (\tr\, W)^2 + 2nu\,\tr\, W + n^2 u^2$.
Moreover,
$\|\tilde{W}\|_F^2 = \|W\|_F^2 + 2u\,\tr\, W + nu^2$.
The first assertion follows.

(ii) Put $W:=\sum_\s\lambda_\s W_\s$. We have
$(\tr W)^2 = \sum_{\s,\t} \lambda_\s\lambda_\t\ (\tr W_\s)(\tr W_\t)$.
By independence and since $W_\s$ are centered we get
$\E\,(\tr W)^2 = \sum_{\s} \lambda_\s^2\ \E\,(\tr W_\s)^2$.
Similarly, we obtain
$\E\,W_{k\ell}^2 = \sum_{\s} \lambda_\s^2\ \E\,((W_\s)_{k\ell})^2$.
Hence $\E\,\|W\|_F^2 = \sum_{\s} \lambda_\s^2\ \|W_\s\|_F^2$.
The second assertion follows.
\end{proof}

Here is a further result stating that stochastic independence follows from invariance.

\begin{lemma}\label{le:invar-indep}
Consider a random $(u,V,W)\in\R\times\R^n\times\Sigma_n$ with joint Gaussian distribution
such that $(u,gV,gWg^T)$ is equivalent to $(u,V,W)$, for all $g\in O(n)$.
(We will call such $(u,V,W)$ $O(n)$-invariant.) Then:
\begin{description}
\item[(i)] $V$ is independent of $u$ and $W$.
\item[(ii)] If $\E\,(u\,\tr\,W)=0$, then $u$ and $W$ are independent.
\end{description}
\end{lemma}

\begin{proof}
$V$ and $W$ are centered by invariance and
we may assume w.l.o.g.\ that $u$ is centered.
(i) Taking $g=-I_n$ we see that $(u,-V,W)$ is equivalent to $(u,V,W)$.
Hence $\E\,(u V_i)=0$ and $\E\,(V_i W_{jk})=0$.

(ii) We argue similarly as in the proof of Proposition~\ref{le:classW}.
Using that $(u,gWg^T)$ is equivalent to $(u,W)$ for
$g=\diag(\e_1,\ldots,\e_n)$ with $\e_i=\pm 1$ we get $\E\,(uW_{ij})=0$
for $i\ne j$. By taking for $g$ a permutation matrix we conclude
$\E\,(uW_{ii})= \E\,(uW_{11})$. By assumption,
$\E\,(uW_{11}) = n^{-1}\E\,(u\,\tr\,W)= 0$.
Hence $u$ is uncorrelated with all~$W_{ij}$.
\end{proof}

The following is a consequence of Lemma~\ref{le:invar-indep} and Gaussian regression.

\begin{corollary}\label{le:crux}
Consider a random $(u,W)\in\R\times\Sigma_n$ with joint Gaussian distribution
such that $(u,gWg^T)$ is equivalent to $(u,W)$, for all $g\in O(n)$.
Let $W_c$ denote the random matrix $W$ conditioned on $u=0$.
Then $W_c$ is $O(n)$-invariant Gaussian and has the following parameter
$$
 \delta(W_c) = \delta(W) - \frac{(\E\,(u\,\tr\,W))^2}{n^2\,\E\,u^2} .
$$
\end{corollary}

\begin{proof}
Consider the random matrix $\tilde{W}:=W - \lambda u I_n$,
where $\lambda:=\frac{\E\,(u\,\tr\,W)}{n\,\E\,u^2}$.
Then $\E\,(u\,\tr\tilde{W})=0$. Moreover,
$(u,\tilde{W})$ is $O(n)$-invariant Gaussian.
According to Lemma~\ref{le:invar-indep},
$u$ and $\tilde{W}$ are independent.
Hence $W_c$ has the same distribution as~$\tilde{W}$.
Moreover, by Lemma~\ref{le:suppositionW},
$$
\delta(\tilde{W}) = \delta(W) + \lambda^2\E\,u^2
 -\frac{2\lambda}{n}\ \E\,(u\,\tr\,W)
 = \delta(W) - \lambda^2\,\E\,u^2
 = \delta(W) - \frac{(\E\,(u\,\tr\,W))^2}{n^2\,\E\,u^2}
$$
as claimed.
\end{proof}

\subsection{Expected determinant of invariant random matrices}\label{se:edet}

This section is devoted to the proof of the following result.
It is stated in~\cite{podk:99}, but  the proof is only vaguely sketched there.

\begin{proposition}\label{pro:expdet}
Suppose $W\in\Sigma_n$ is $O(n)$-invariant and Gaussian, $n>1$.
Then we have
$$
 \E\, \det(I_n + W) = \E\, \big( 1 + \sqrt{\delta(W)}\, X\big)^n,
$$
where $X\sim N(0,1)$ is a standard normal random variable.
(Note that $\delta(W)$ may be negative.)
\end{proposition}

To prepare for the proof we denote by ($k\in\N$)
\begin{equation}\label{eq:def-gamma}
\gamma_k =\E\, |X|^k = \frac2{\sqrt{2\pi}}\int_0^\infty y^k e^{-\frac{y^2}{2}}\, dy
         =\frac1{\sqrt{\pi}}\, 2^{k/2}\, \Gamma\left(\frac{k+1}{2}\right)
\end{equation}
the $k$-th absolute moment of a standard normal random variable $X\sim N(0,1)$.
In particular,
\begin{equation}\label{eq:gamma}
\gamma_{2m} = 1\cdot 3\cdot 5\cdots(2m-1) = (2m)!/(2^m m!),\quad
\gamma_k\Oh_k = 2(2\pi)^{k/2}.
\end{equation}

The following lemma about the higher moments of Gaussian vectors
is well known, cf.~\cite[p.~108]{adl:81}.

\begin{lemma}\label{le:moments}
Let $(Y_1,\ldots,Y_n)$ be a centered Gaussian vector.
Then
$\E (Y_1\cdots Y_n) =0$
if $n$ is odd.
Otherwise, if $n=2m$,
$$
 \E (Y_1\cdots Y_{2m}) = \sum \E (Y_{i_1}Y_{i_2}) \cdots \E (Y_{i_{2m-1}}Y_{i_{2m}}) ,
$$
where the sum is over all $\gamma_{2m}=(2m)!/(2^m m!)$ different ways of grouping
$Y_1,\ldots,Y_{2m}$ into $m$ pairs.
\end{lemma}

Clearly, the lemma implies that $\E\,\det W = 0$ for any centered Gaussian matrix
$W\in\R^{n\times n}$ if $n$ is odd.

\begin{corollary}\label{cor:EdetT}
The random matrix $T\in\Sigma_{n}$ of Proposition~\ref{le:classW} satisfies
$\E\,\det T = (-1)^m \gamma_{2m}$ if $n=2m$ and $\E\,\det T =0$ otherwise.
\end{corollary}

\begin{proof}
We have
$\E\,\det T = \sum_{\pi}\, \mathrm{sgn}(\pi)\,\E (T_{1,\pi(1)}\cdots T_{2m,\pi(2m)})$.
Lemma~\ref{le:moments} and the fact that the entries of $T$
are independent imply that only the products of transpositions of the form
$\pi = (i_1 j_1)\cdots (i_m j_m)$ with
$\{i_1,j_1,\ldots,i_m,j_m\} =\{1,2,\ldots,2m\}$
contribute to the sum.
The contribution of each such $\pi$ is $(-1)^m$ and there are
$\gamma_{2m}$ such $\pi$.
\end{proof}

\proofof{Proposition~\ref{pro:expdet}}
According to Proposition~\ref{le:classW} we may assume
$W=rZ I_n + sT$. Then $\delta(W)=r^2-s^2$.
By expanding the determinant it is easy to see that
$$
 \det (I_n +W) = \det ((1+rZ)I_n + sT)
  = \sum_K (1 + rZ)^{n-|K|}\, s^{|K|}\, \det T_{K,K},
$$
where the sum is over all subsets $K$ of $\{1,2,\ldots,n\}$
and $T_{K,K}$ denotes the principal submatrix obtained by
selecting the rows and columns in $K$.
Using Corollary~\ref{cor:EdetT} and taking into account the
independence of $Z$ and $T$ we obtain
$$
 \E\,\det (I_n +W) =
 \sum_{j=0}^{\lfloor n/2\rfloor} {n\choose 2j} \E\,(1 + rZ)^{n-2j}\,\, s^{2j} (-1)^j \gamma_{2j} .
$$
Suppose $Y\sim N(0,1)$ is independent of $Z$. Using $i=\sqrt{-1}$
we have $(-1)^j\gamma_{2j} = \E (iY)^{2j}$ and hence
\begin{eqnarray*}
 \det (I_n +W)
 &=& \sum_{j=0}^{\lfloor n/2\rfloor} {n\choose 2j} \E\,(1 + rZ)^{n-2j}\,\E\, (isY)^{2j} \\
 &=& \E (1 + rZ + isY)^n \\
 & =& \sum_{j=0}^{\lfloor n/2\rfloor} {n\choose 2j} \E\,(rZ+ isY)^{2j} .
\end{eqnarray*}
On the other hand, using that
${2j\choose 2k}\,\gamma_{2k}\gamma_{2j-2k}={j\choose k}\,\gamma_{2j}$ for $0\le k\le j$, we get
\begin{eqnarray*}
 \E\,(rZ+ isY)^{2j} &=& \sum_{k=0}^{j}{2j\choose 2k}\, r^{2k}\gamma_{2k}\,
                        (-1)^{j-k}\,s^{2j-2k}\, \gamma_{2j-2k} \\
                    &=& \gamma_{2j} \sum_{k=0}^{j}{j\choose k}\, r^{2k}\,(-1)^{j-k}\,s^{2j-2k}
                             = \gamma_{2j} (r^2-s^2)^j .
\end{eqnarray*}
We conclude that
$$
 \det (I_n +W) = \sum_{j=0}^{\lfloor n/2\rfloor} {n\choose 2j} \E\,(\sqrt{r^2-s^2} Z)^{2j}
               = \E\,( 1 + \sqrt{r^2-s^2} Z)^n ,
$$
which finishes the proof.
\proofend

\subsection{Two auxiliary results on random matrices}\label{se:sol-rand-le}
Let $A\in\R^{s\times n}$ be a matrix of rank~$s\le n$.
We denote by $A_{m}\in\R^n$ the $m$-th row of $A$ and by
$A_m^\perp$ the orthogonal projection of $A_m$ onto the space spanned
by $A_1,\ldots,A_{m-1}$.
Let $\vol(A_1,\ldots,A_m)$ denote the volume of the parallelepiped in $\R^n$
spanned by $A_1,\ldots,A_m$.
Clearly,  $\vol(A_1,\ldots,A_m) =  \vol(A_1,\ldots,A_{m-1})\, \|A_m^\perp\|$.

\begin{lemma}\label{le:techcalc}
Suppose that $A\in\R^{s\times n}$ is a random matrix with independent
standard Gaussian entries. Then we have for $1\le m\le s$
and $0\le j\le \lfloor\frac{n-m}{2}\rfloor$
$$
 \E_A\Big(\vol(A_1,\ldots,A_{m-1}) \cdot \frac1{\|A_m^\perp\|^{2j-1}}\Big)
  = \frac{\Oh_{n-m}} {2 (2\pi)^{\frac{n-m}{2}}}\ \frac{\gamma_n \gamma_{n-m+1-2j}}{\gamma_{n-m+1}}.
$$
\end{lemma}

\begin{proof}
The invariance under rotations of the standard normal distribution shows that,
conditioning on $A_1,\ldots,A_{m-1}$, the random variable $A_m^\perp$
has a standard normal distribution in the orthogonal complement of $A_1,\ldots,A_{m-1}$
in $\R^n$, which is of dimension $n-m+1$ with probability one.
Hence
$\E\,\big(\|A_m^\perp\| \big/ A_1,\ldots,A_{m-1}\big) = K_{n-m+1}$,
where we have written $K_d:=\E\,\|X\|$ for a standard normal $X$ in $\R^d$.
An elementary computation gives (cf.~\cite{azw:05})
\begin{equation*}\label{eq:KL-formula}
 K_d = \sqrt{2}\, \frac{\Gamma((d+1)/2)}{\Gamma(d/2)},\quad
 K_1 K_2 \cdots K_{m} = \frac1{\sqrt{\pi}}\, 2^{\frac{m}{2}}\ \Gamma\left(\frac{m+1}{2}\right) = \gamma_{m} .
\end{equation*}
Put $\vol_m(A):=\vol(A_1,\ldots,A_m) $. Then
$\vol_m(A) =\vol_{m-1}(A)\, \|A_m^\perp\|$.
Hence
$$
 \E\,\vol_m(A) = \E\left(\vol_{m-1}(A)\ \E\,\left(\|A_m^\perp\|  \Big/ A_1,\ldots,A_{m-1} \right) \right)
  = \E\,\vol_{m-1}(A)\ K_{n-m+1} ,
$$
which implies
\begin{equation}\label{eq:volA}
 \E\,\vol_m(A) = K_{n-m+1} K_{n-m+2}\cdots K_n = \frac{\gamma_n}{\gamma_{n-m}}.
\end{equation}
Using polar coordinates,
we get for almost all values of $A_1,\ldots,A_{m-1}$
\begin{eqnarray*}
 \E_{A_m} \left( \|A_m^\perp\|^{1-2j} \Big/ A_1,\ldots,A_{m-1}\right)
 &=& \frac{\Oh_{n-m}}{ (2\pi)^{\frac{n-m+1}{2}}} \int_0^\infty r^{n-m-2j+1} e^{-\frac{r^2}{2}}\, dr \\
 &=& \frac12\, \frac{\Oh_{n-m}}{ (2\pi)^{\frac{n-m}{2}}}\, \gamma_{n-m-2j+1} .
\end{eqnarray*}
The claim follows by combining this with (\ref{eq:volA}).
\end{proof}

The Moore-Penrose inverse of a matrix $A\in\R^{s\times n}$ of rank~$s$ is defined by
$A^\dagger := A^T (A A^T)^{-1} \in \R^{n\times s}$,
cf.~\cite{bell:97}. 
It is characterized by the following properties:
$AA^\dagger=I_s$ and $A^\dagger A$ is the orthogonal projection onto
the orthogonal complement $(\ker A)^\perp$ of the kernel of~$A$.
Note that $(\ker A)^\perp$ is generated by the rows of $A$.
Let us denote by $\tilde{A}$ the restriction of the linear map $A\colon\R^n\to\R^s$
to $(\ker A)^\perp$.
Then the linear map $A^\dagger\colon\R^s\to\R^n$ is the inverse of $\tilde{A}$.
We note that
\begin{equation}\label{eq:volvol}
|\det\tilde{A}| = \vol(A_1,\ldots,A_s) =\sqrt{\det (A A^T)}.
\end{equation}
(Proof:  $\tilde{A}^Te_i = A_i$, hence
$|\det\tilde{A}^T| = \vol(A_1,\ldots,A_s)$. It is well known that
the latter equals  $\sqrt{\det A A^T}$.)

We denote by $(N_1,\ldots,N_s)$ the orthonormal basis of $(\ker A)^\perp$
obtained from the rows of $A$ by Gram-Schmidt orthogonalization.
This defines  the orthogonal map
\begin{equation}\label{eq:QA}
Q_A\colon \R^s\to (\ker A)^\perp,\ u\mapsto \sum_{\s=1}^s u_\s N_\s .
\end{equation}
By composing $Q_A$ with the adjoint
$(A^\dagger)^T$ of the Moore-Penrose inverse~$A^\dagger$
we obtain the linear endomorphism
$(A^\dagger)^T Q_A \colon \R^s\to\R^s$.

The following is the main result of this subsection.
It will be used in the proof of Theorem~\ref{th:main}.

\begin{proposition}\label{le:heureka}
Suppose that $A\in\R^{s\times n}$ is a random matrix with independent
standard Gaussian entries. Let $u\in S^{s-1}$ be a uniformly distributed
random unit vector which is independent of $A$.
Then the random variable
$$
 \Big ( \sqrt{\det A A^T},   (A^\dagger)^T Q_A (u) \Big)
$$
with values in $\R\times\R^s$ has the same distribution as
$$
 \Big( \vol(A_1,\ldots,A _{s-1})\cdot \| A_s^\perp\|, \frac1{\| A_s^\perp\|} \cdot w \Big)
$$
where $w$ is uniformly distributed in the sphere $S^{s-1}$ and
independent of $A$.
\end{proposition}

We give some preparations for the proof.
Let $A\in\R^{s\times n}$ be of rank~$s$.
Pick $y\in (\ker A)^\perp$ and define
$A^{(y)}$ as the restriction of $A$ to the orthogonal complement of
$\R y \oplus \ker A$. This can be described by the following
commutative diagram (where upgoing arrows denote inclusions):
$$
\begin{array}{ccc}
 (\ker A)^\perp & \stackrel{\tilde{A}}{\longrightarrow} & \R^s \\
        \uparrow & &  \uparrow \\
 (\R y \oplus \ker A)^\perp &\stackrel{A^{(y)}}{\longrightarrow} & \im A^{(y)}
\end{array}.
$$
Suppose now $\|y\|=1$ and
let $(Ay)^\perp$ denote the orthogonal projection of $Ay$ onto
the image of $A^{(y)}$. Then it is easy to see that
\begin{equation}\label{eq:volly}
 |\det \tilde{A}| = |\det A^{(y)}|\cdot \|(Ay)^\perp\| .
\end{equation}

\begin{lemma}\label{le:dagger}
Let $A\in\R^{s\times n}$ be of rank~$s\le n$ and $y\in(\ker A)^\perp$
be a unit vector. Then the length of the vector $v=(A^\dagger)^T y\in\R^s$
can be described as
$$
 \|v \| = \left|\frac{\det A^{(y)}}{\det\tilde{A}}\right| = \frac1{\|(Ay)^\perp \|} .
$$
\end{lemma}

\begin{proof}
We have $v^T A = y^T A^\dagger A = y^T$, hence $A^Tv=y$.
Thus $v$ can be interpreted as the unique solution of the system of linear equations
$$
 v_1 A_1 + \cdots + v_s A_s = y .
$$
Cramer's rule implies that
$$
 | v_i| = \frac{\vol(A_1,\ldots,A_{i-1},y,A_{i+1},\ldots,A_s)}{\vol(A_1,\ldots,A_s)} .
$$
Let $\overline{A}_i$ denote the projection of $A_i$ onto the orthogonal complement
of $\R y$. Then
$$
 \vol(A_1,\ldots,A_{i-1},y,A_{i+1},\ldots,A_s) =
 \vol(\overline{A}_1,\ldots,\overline{A}_{i-1},\overline{A}_{i+1},\ldots,\overline{A}_s) .
$$
It is now sufficient to prove that (cf.~(\ref{eq:volvol}))
\begin{equation}\label{eq:binet}
 \sum_{i=1}^s \vol(\overline{A}_1,\ldots,\overline{A}_{i-1},\overline{A}_{i+1},\ldots,\overline{A}_s) ^2
  = (\det A^{(y)})^2 .
\end{equation}
By an orthogonal transformation, we may assume without loss of generality that
$\ker A=0\times \R^{n-s}$ and $y=(0,\ldots,0,1,0,\ldots,0) $ (with $1$ at position~$s$).
Write $A=(a_{ij})$. Then $a_{ij}=0$ for $j>s$ and
$\tilde{A}$ is given by the matrix $(a_{ij})_{i,j\le s}$ with respect to the canonical bases.
Moreover, $A^{(y)}$ is given by the matrix $M:=(a_{ij})_{i\le s,j <s}$ in $\R^{s\times (s-1)}$.
The relation of Binet-Cauchy~\cite{bell:97} states that
$$
 \sum_{i=1}^s (\det M_i )^2 = \det (M^T M) ,
$$
where $M_i$ stands for the square matrix obtained from $M$ by deleting the
$i$-th row. This is exactly the asserted Equation~(\ref{eq:binet}).
\end{proof}

\proofof{Proposition~\ref{le:heureka}}
Consider the following random process. First choose $A\in\R^{s\times n}$
at random with independent standard Gaussian entries.
Then pick a unit vector $y\in(\ker A)^\perp$ uniformly at random.
Then the resulting random variable
$\big ( \sqrt{\det A A^T},   (A^\dagger)^T (y)\big)$
is equivalent to
$\big ( \sqrt{\det A A^T},   (A^\dagger)^T Q_A (u) \big)$.

Equations~(\ref{eq:volvol}) and (\ref{eq:volly}) give that
$$
 \sqrt{\det (A A^T)} = |\det\tilde{A}| = |\det A^{(y)}| \cdot \|(Ay)^\perp\| .
$$
Lemma~\ref{le:dagger} implies that
$$
 (A^\dagger)^T (y)  = \frac1{\|(Ay)^\perp \|} \cdot w
$$
with a random unit vector $w\in S^{s-1}$.
It is clear that
$(A^\dagger)^T (y)$ is $O(s)$-invariant.
According to Lemma~\ref{le:univect},
this implies that $w$ is uniformly distributed and independent of
$\|(Ay)^\perp\|$.

As in the proof of Lemma~\ref{le:techcalc}, we see that
$\|A_s^\perp\|$ conditioned on $A_1,\ldots,A_{s-1}$
is equivalent to $\|X\|$, where $X$~is standard normal in $\R^{n-s+1}$.
It is not hard to see that
$\|(Ay)^\perp\| $ has the same distribution as $\|X\|$.
Similarly, one shows that $|\det A^{(y)}| $ is equivalent to
$\vol(A_1,\ldots,A_{s-1})$.
(Note that if $g_1,\ldots,g_{s-1},y$ is an orthonormal basis of $(\ker A)^\perp$,
then $|\det A^{(y)}| = \vol(Ag_1,\ldots,Ag_{s-1})$.)
Finally,  $\|(Ay)^\perp\| $ is independent of $|\det A^{(y)}| $.
\proofend

\section{Invariant random polynomials}\label{se:invarpoly}

\subsection{Classification}\label{se:classific}

We briefly describe the classification of invariant Gaussian polynomials
and refer to Kostlan~\cite{kost:93,kost:02} for more details.
This classification is just for illustration and will not be needed for
the proof of Theorem~\ref{th:main}.

Recall that $H_{d,n}$ denotes the vector space of homogeneous real polynomials of degree~$d$
in the variables $X_0,\ldots,X_n$.
The orthogonal group $O(n+1)$ operates
on $H_{d,n}$ (from the left) in the natural way:
for $f\in H_{d,n}$ and $g\in O(n+1)$ we set $(gf)(X):=f(g^{-1}X)$,
where $X=(X_0,\ldots,X_n)^\top$.

\begin{definition}\label{def:invar-rand-poly}
A random polynomial $f\in H_{d,n}$ is called
{\em $O(n+1)$-invariant}  (invariant for short)
iff $gf$ is equivalent to $f$ for all all $g\in O(n+1)$.
\end{definition}

\begin{remark}
Any invariant random polynomial $f\in H_{d,n}$ is centered if $d$ is odd.
Moreover, if $d$ is even, it is easy to reduce to the case where $f$ is centered,
cf.~\cite[\S5.3]{kost:02}.
For convenience we will therefore additionally require that $f$ is centered.
\end{remark}

\begin{example}\label{ex:kostlan}
The most natural example of an invariant random polynomial in $H_{d,n}$
is obtained as follows.
Write $f\in H_{d,n}$ in the form
$f=\sum_\a f_\a X^\a$, where the sum is over all multiindices
$\a\in\N^{n+1}$ such that $|\a|:=\sum_i \a_i = d$.
Assume that the coefficients $f_\a$ are independent with
centered Gaussian distribution and variance
$\Var(f_\a) = {d\choose \a} := \frac{d!}{\a_0!\cdots\a_n!}$.
The covariance function
$\R^{n+1}\times\R^{n+1}\to\R, (x,y)\mapsto\E\,(f(x)f(y))$
satisfies
$$
\E\,(f(x)f(y)) = \sum_{\a,\b}\E(f_\a f_\b) x^\a y^\b
  =\sum_{\a} {d\choose\a}x^\a y^\a
  = (\sum_{i=0}^n x_i y_i)^d = \langle x,y\rangle^d .
$$
and it is thus invariant under the action of $O(n+1)$. It follows that
$f$~is an invariant Gaussian random polynomial.
We will say that $f$ is {\em Kostlan distributed}~\cite{kost:93}.
Up to a scalar, this random polynomial can be characterized by requiring
that $f\in H_{d,n}$ is invariant and centered
and has stochastically independent coefficients $f_\a$
(cf.~\cite[Theorem 4.5]{kost:93}). It is interesting that
the normal distribution is enforced by the above requirements only.
\end{example}

It is possible to characterize all invariant centered Gaussian random
polynomials of $H_{d,n}$.
For doing so, it is helpful to start with some general observations.

A centered Gaussian distribution on $V=\R^n$ is characterized
by its covariance matrix. In coordinate free language, this corresponds
to the choice of an inner product on the dual space $V^\ast$ defined by
$$
V^\ast\times V^\ast \to \R,\ (\lambda,\mu)\mapsto \E_{\,v}\,( \lambda(v)\mu(v))
$$
Consider now the general situation of a compact Lie group $G$ operating on
a real vector space~$V$ (in our situation $G=O(n+1)$).
It is easy to see that a centered Gaussian distribution on $V$
is $G$-invariant iff the corresponding inner product is $G$-invariant, that is,
$\langle gu,gv \rangle = \langle u,v \rangle$ for all $g\in G$ and all $u,v\in V$.

From the representation theory of compact groups it is known that
\begin{enumerate}
\item
On an irreducible $G$-module $W$ there is a (up to a positive scalar)
unique $G$-invariant inner product.

\item Any two nonisomorphic submodules of $V$ are orthogonal with respect
to a $G$-invariant inner product of $V$.
\end{enumerate}
(The first statement follows from Schur's lemma; for the second see
\cite[\S27, p.~29]{hero:70}.)

Suppose that $V$ splits into a direct sum $V=\oplus_{i=1}^m W_i$ of pairwise nonisomorphic
irreducible submodules $W_i$. Choose an invariant inner product $\langle\ ,\ \rangle_i$
on each $W_i$. Then, according to the above facts,
the invariant inner products on $V$ are of the form
$$
 \mbox{$V\times V \to\R, (\oplus u_i, \oplus_i v_i) \mapsto
  \sum_{i=1}^m c_i \langle u_i , v_i \rangle_i$},
$$
parameterized by $c_1,\ldots,c_m>0$.

We briefly describe the decomposition of the $O(n+1)$-module $H_{d,n}$ into irreducible
submodules, cf.~\cite[\S5.2.3]{gowa:98}.
Write $r^2 := \sum_{i=0}^n X_i^2$ and consider the equivariant linear map
$\varphi\colon H_{d,n} \to H_{d,n}, f\mapsto r^2 \Delta f$ arising from
the Laplace operator $\Delta$ on $\R^{n+1}$.
It is known that the kernel $\mathcal{H}_{d,n}$ is irreducible.
The elements of $\mathcal{H}_{d,n}$ are called {\em harmonic polynomials} of degree~$d$.
The map $\varphi$ has the following eigenspace decomposition
\begin{equation}\label{eq:harmpd}
 H_{d,n} = \bigoplus_{i=0}^{\lfloor d/2\rfloor} r^2 \mathcal{H}_{d-2i,n} ,
\end{equation}
which is thus a decomposition of $H_{d,n}$ into nonisomorphic irreducible submodules
($\mathcal{H}_{d-2i,n} $ corresponds to the eigenvalue $2i(n+2d-2i-1)$).
We conclude that the invariant inner products on $H_{d,n}$ can be
parameterized by $\lfloor d/2\rfloor$ positive numbers.
We refrain from explicitly describing these inner products and refer instead to
\cite{kost:02} for details.
To summarize, we see that the invariant centered Gaussian random
polynomials of $H_{d,n}$ can be parameterized by $\lfloor d/2\rfloor$ positive numbers.

\subsection{Parameter of invariant random polynomials}\label{se:par-randpoly}

A polynomial $f\in H_{d,n}$ defines a differentiable map
$f\colon S^n\to\R$.
We denote by $Df(x)\colon T_xS^n\to\R$ and
$D^2f(x)\colon T_x S^n\times T_x S^n\to\R$
the first and second order derivative of $f$ at $x\in S^n$.
They are characterized by
\begin{equation}\label{eq:taylor}
 f\bigg(\frac{x + \lambda V}{\|x + \lambda V\|}\bigg) =
 f(x) + \lambda Df(x)(V) + \frac{\lambda^2}{2} D^2f(x)(V,V) + o(\lambda^2)
\end{equation}
for $V\in T_xS^n$ and $\lambda\in\R$ going to zero.
At the point $q=(1,0\ldots,0)$ we can identify
the tangent space $T_qS^n = 0\times\R^n$ with $\R^n$.
Clearly, $Df(q)$ is determined by the gradient vector
$(\partial_1 f(q),\ldots,\partial_n f(q))$.

\begin{lemma}\label{le:2nd-deriv}
$D^2f(q)$ is given by the matrix
$(\partial^2_{k\ell} f(q))_{1\le k,\ell\le n} - d\, f(q)\, I_n$.
\end{lemma}

\begin{proof}
Consider
$g(t):=f((1+t^2)^{-1/2}(1,t V))= (1+t^2)^{-d/2} f(1,tV)$
for $V\in\R^n$ of length~$1$.
By definition of the second order derivative we have $D^2f(q)(V,V) = g''(0)$,
cf.~(\ref{eq:taylor}).
A calculation yields
$g''(0) = \sum_{k\ell} \partial^2_{k\ell}f(q) V_k V_\ell -d\, f(q)$.
\end{proof}

The following definition is from~\cite{podk:99}.

\begin{definition}\label{def:invar-poly}
The {\em parameter} $\d(f)$ of  an invariant random polynomial $f\in H_{d,n}$ is defined as
$$
  \delta(f):=\frac{\|x\|^2\, \E\,\|Df(x)\|^2}{n\E\,f(x)^2}
$$
where $x\in\R^{n+1}\setminus \{0\}$
(this is independent of $x$ by invariance and homogeneity).
\end{definition}

\begin{remark}\label{re:par-f}
\begin{enumerate}
\item If $f\in H_{d,n}$ is invariant then $Df(q)$ is
Gaussian with covariance matrix $\delta(f) \, \E\,f(q)^2\, I_n$.
Thus $\delta(f)=\frac{\E\, (\partial_{k} f(q))^2}{\E\,f(q)^2}$
for any $1\le k\le n$.

\item The parameter of the Kostlan distribution equals the degree~$d$.

\end{enumerate}
\end{remark}

\begin{proof}
The first statement follows from Lemma~\ref{le:invar-vect}.
For the second just note that $f(q)=f_{(d,0,\ldots,0)}$ and
$\partial_{X_1}f(q)=f_{(d-1,1,0,\ldots,0)}$.
\end{proof}


\begin{lemma}\label{le:restrict-invar-distr}
Suppose $f\in H_{d,n}$ is $O(n+1)$-invariant and $n'\le n$.
Then the restriction $f'\in H_{d,n'}$ of $f$ to $\R^{n'+1}$ is $O(n'+1)$-invariant
and has the same parameter, i.e., $\delta(f')=\delta(f)$.
Moreover, if $f$ is Kostlan distributed, then so is~$f'$.
\end{lemma}

\begin{proof}
It is clear that the distribution of $f'$ is invariant with
respect to the action of $O(n'+1)$. In order to see that the
parameter remains the same just note that restricting $f$ to $\R^{n'+1}$
means substituting the variables $X_{n'+1},\ldots,X_n$ by $0$.
The assertion about the Kostlan distribution is obvious.
\end{proof}

We show now that the first and second order derivatives of $f$ inherit
the invariance property from $f$.

\begin{lemma}\label{le:dervq-invar}
Let $f\in H_{d,n}$ be an invariant Gaussian random polynomial.
Then $(f(q),Df(q),D^2f(q))\in\R\times\R^n\times\Sigma_n$
is $O(n)$-invariant in the sense of Lemma~\ref{le:invar-indep}.
In particular, $Df(q)$ is independent of $f(q)$ and $D^2f(q)$.
\end{lemma}

\begin{proof}
Consider for fixed $g\in O(n+1)$ the transformed polynomial
$h:=gf$, that is, $h(x)= f(g^T x)$.
By assumption, $h$ has the same distribution as $f$.
This implies that the random vector
$(h,\partial_k h,\partial^2_{k\ell}h)$ is equivalent to
$(f,\partial_k f,\partial^2_{k\ell}f)$.
(This can be shown by expanding $f$ and $g$ in a fixed basis of
$H_{d,n}$ with random real coefficients.)
We conclude that
$(h(q),Dh(q),D^2h(q)) = (f(q),gDf(q),gD^2f(q)g^T)$
is equivalent to $(f(q),Df(q),D^2f(q))$.
\end{proof}

The following proposition from~\cite{podk:99} says that the parameter of $f\in H_{d,n}$
equals the parameter of $D^2f(q)$, up to a scaling factor.
Since the proof in \cite{podk:99} is incomplete,
we provide a different proof in the appendix.
The assumption of a Gaussian random polynomial is only
made for simplifying the statement and could be replaced by
suitable regularity conditions.

\begin{proposition}\label{pro:delta=}
Suppose that $f\in H_{d,n}$ is an invariant Gaussian random polynomial. Then:
\begin{enumerate}
\item[(i)] $\E\,(f(q)\,\tr\, D^2f(q)) = - n\delta(f)\,\E\,f(q)^2$.

\item[(ii)] $\delta (D^2f(q)) = \delta(f)\ \E\, f(q)^2$.
\end{enumerate}
\end{proposition}

The next corollary will be crucial in the proof of Theorem~\ref{th:main}.

\begin{corollary}\label{cor:delta-cond-appl}
Suppose that $f\in H_{d,n}$ is an invariant Gaussian random polynomial.
Let $W_c$ be the random matrix $D^2f(q)$ conditioned on $f(q)=0$. Then we have
$\delta(W_c) = \delta(f)(1-\delta(f))\,\E\,f(q)^2$.
\end{corollary}

\begin{proof}
Put $u=f(q)$ and $W=D^2f(q)$.
Lemma~\ref{le:dervq-invar} implies that $(u,W)$ is
$O(n)$-invariant.
Proposition~\ref{pro:delta=} yields $\delta(W)=\delta(f)\ \E\,u^2$ and
$\E\,(u\,\tr\, W) = - n\delta(f)\,\E\,u^2$.
The assertion follows now from Corollary~\ref{le:crux}
\end{proof}

\section{Random real projective varieties}\label{se:randvar}

We give here the proof of the main Theorem~\ref{th:main}.
Starting from Weyl's tube formula~(\ref{eq:weyl})
we present in \S\ref{se:rice}
a version of a ``Rice formula'' for curvature coefficients.
We then proceed with a probabilistic analysis of that formula,
making heavily use of the invariance under the orthogonal group.
In order to do so, we need all the auxilary material on invariant random
vectors, matrices, and polynomials that was collected in
\S\ref{se:vectpolymat}--\S\ref{se:invarpoly},
except~\S\ref{se:classific}.

We remark that the kinematic formula of integral geometry
(Theorem~\ref{th:kinform-S}) would allow to reduce to the considerably simpler case of one equation.
However, in view of a further development of the theory (higher moments),
we will not use the kinematic formula here but instead give
a self-contained probabilistic proof.

\subsection{A Rice formula for expected curvature coefficients}\label{se:rice}

In a first step we are going to derive a somewhat more explicit form
of Weyl's tube formula~(\ref{eq:weyl}) for the zero set of homogeneous
polynomials in $S^n$.

Let $f_1,\ldots,f_s\in\R[X_0,\ldots,X_n]$ be homogeneous real polynomials
of the degrees $d_1,\ldots,d_s$ ($1\le s\le n$).
They define a differentiable map $f\colon S^n\to\R^s$.
For a point $x\in S^n$ we denote by
$Df(x)\colon T_xS^n\to\R^s$ and
$D^2f(x)\colon T_x S^n\times T_x S^n\to\R^s$
the first and second order derivative of $f$ at $x$.

In the following we assume that
$x\in S^n$, $f(x)=0$, and $\mathrm{rank}\, Df(x) =s$.
Then, locally at~$x$, the zero set $Z=\mZ(f)$ is a smooth Riemannian
submanifold of $S^n$ of dimension $n-s$.
The kernel of $Df(x)$ equals the tangent space $T_xZ$ of $Z$ at~$x$.
We denote the inverse of the restriction of $Df(x)$ to the orthogonal
complement $(T_xZ)^\perp$ of $T_xZ$ by
$Df(x)^\dagger\colon\R^s\to (T_xZ)^\perp$
(Moore-Penrose inverse, compare~\S\ref{se:sol-rand-le}).

The following, certainly well known lemma, expresses the second fundamental form
$\zwff_Z(x)$ of $Z$ at $x$ (cf.~\S\ref{se:weyltube})
in terms of the {\em Hessian} $Hf(x)\colon T_xZ\times T_xZ\to\R^s$, which we define as
the restriction of $D^2f(x)$ to $T_xZ\times T_xZ$.
Since we could not find an appropriate reference,
we have included a proof in the appendix.

\begin{lemma}\label{le:2ff-implicit}
Under the above assumptions we have for $V,W\in T_xZ$ and $\nu\in (T_xZ)^\perp$
$$
 \langle L_Z(x,\nu)(V),W\rangle
 = \zwff_Z(x)(V,W,\nu) = -\langle \nu, Df(x)^\dagger Hf(x)(V,W)\rangle.
$$
\end{lemma}

The derivative $Df_\s(x)\colon T_xS^n\to\R$ can be identified with a vector in
$T_xS^n$ via the inner product on $T_xS^n$. Suppose $(N_1,\ldots,N_s)$ is the
orthonormal basis of $(T_xZ)^\perp$ obtained from
$(Df_1(x),\ldots,Df_s(x))$ by Gram-Schmidt orthogonalization. We use the
orthogonal map
$Qf(x)\colon\R^s\to (T_xZ)^\perp,\ u\mapsto \sum_{\s=1}^s u_\s N_\s$
to describe unit normal vectors in $(T_xZ)^\perp$ by coordinates.
We thus define a Weingarten map
$$
  Lf(x,u) := L_Z(x,Qf(x)(u))
$$
of $Z$ at $x$ in direction parameterized by $u\in S^{s-1}$ (compare \S\ref{se:weyltube}).
According to Lemma~\ref{le:2ff-implicit},
the Weingarten map is explicitly characterized by
\begin{eqnarray}\notag
 \langle Lf(x,u)(V),W\rangle
                             &=& - \langle Qf(x)(u), Df(x)^\dagger Hf(x)(V,W)\rangle\\ \label{eq:defL}
                             &=& -\langle (Df(x)^\dagger)^T Qf(x)(u), Hf(x)(V,W)\rangle,
\end{eqnarray}
where $(Df(x)^\dagger)^T\colon (T_xZ)^\perp\to\R^s$ denotes the adjoint map of
$Df(x)^\dagger\colon\R^s\to (T_xZ)^\perp$.

We suppose now that $\mathrm{rank}\, Df(x) =s$ for all $x\in Z=\mZ(f)$, i.e.,
the hypersurfaces $\mZ(f_i)$ intersect transversally.
By Sard's lemma, this is the case for almost all $f$.
Then $Z$ is either empty or a compact smooth submanifold of $S^n$ of dimension~$n-s$ 
and we assume the latter. 

It will be convenient to introduce the following function associated with~$f$
$$
 g_f\colon S^n\times \R^s \to\R,\ g_f(x,t)
  := \frac{\det(\id - \|t\| Lf(x,t/\|t\|))}{(1+ \|t\|^2)^{(n+1)/2}} .
$$
By the transformation theorem, Weyl's formula~(\ref{eq:weyl}) for the volume of tubes
can be concisely rewritten as
$\vol(T(Z,\a)) = \int_{Z\times B_a} g_f\, d (Z\times B_a)$, 
where $a=\tan\a$ is sufficiently small
and $B_a:=\{t\in\R^s\mid \|t\| \le a \}$
is the ball of radius $a$ in $\R^s$.
Combining this with (\ref{eq:voltube})  we obtain
\begin{equation}\label{eq:weylexpand}
   \int_{Z\times B_a} g_f\, d (Z\times B_a)
   = \sum_{j=0}^{\lfloor\frac{n-s}{2}\rfloor} K_{s+2j}(Z)\, J_{n,s+2j}(\a) .
\end{equation}
This expansion is valid for all $0<\a<\pi/2$ since the functions on both sides
of~(\ref{eq:weylexpand}) are analytic.

We can now state the announced Rice formula for expected curvatures, 
which will allow to determine the expectations $\E_f K_{s+2j}(\mZ(f))$
of the curvature coefficients.

\begin{theorem}\label{th:rice}
Suppose that $f=(f_1,\ldots,f_s)\in H_{d_1,n}\times\cdots\times H_{d_s,n}$
is a Gaussian random vector.
Hence $f(x)\in\R^s$ is a Gaussian random vector for any $x\in S^n$,
and we shall denote its density function by
$p_{f(x)}\colon\R^s\to\R$.
Let $u\in S^{s-1}$ be a random unit vector uniformly distributed in the sphere
and independent of $f$.
Define the function
$\psi\colon [0,\infty)\times S^n\to\R$
by the following conditional expectation for
$(r,x)\in [0,\infty)\times S^n$
\begin{equation*}
 \psi(r,x) :=  p_{f(x)}(0)\,\E_{f,u}\left(\sqrt{\det(Df(x)Df(x)^T)}
  \det(\id - r Lf(x,u))\Big/ f(x)=0 \right).
\end{equation*}
By taking a spherical average we define the function
$$
 \Psi\colon [0,\infty)\to\R,\
 \Psi(r):= \frac1{\Oh_n}\int_{x\in S^n} \psi(r,x)\, dS^n.
$$
Then we have for $0<\a<\pi/2$ and $a=\tan\a$
\begin{gather*}
 \sum_{j=0}^{\lfloor\frac{n-s}{2}\rfloor} \E\, K_{s+2j}(\mZ(f))\, J_{n,s+2j}(\a)
 = \Oh_n\Oh_{s-1} \int_{0}^a \frac{r^{s-1}\Psi(r)}{(1+r^2)^{(n+1)/2}}\, dr .
\end{gather*}
\end{theorem}

\begin{proof}
The proof uses similar ideas 
as in \cite[\S 5.1]{adl:81} and \cite[Theorem~1]{azw:05}.

Fix $a>0$ and consider for 
$f_i\in H_{d_i,n}$, $1\le i\le s$, 
the corresponding map $f\colon S^n\to\R^s$. 
The fibre integral
$$
 G_f(y) := \int_{f^{-1}(y)\times B_a} g_f\, d (f^{-1}(y)\times B_a)
$$
is well defined for regular values $y\in\R^s$. 
We thus need to determine (cf.~(\ref{eq:weylexpand}))
\begin{equation}\label{eq:globy}
\sum_{j=0}^{\lfloor\frac{n-s}{2}\rfloor} \E K_{s+2j}(\mZ(f))\, J_{n,s+2j}(\a) 
 = \E_f (G_f(0)) .
\end{equation}

We will apply the coarea formula (or Fubini's theorem for 
Riemannian manifolds).  
Recall that the {\em normal Jacobian} $\NJ_f(x) :=\sqrt{\det(Df(x) Df(x)^T)}$
of $f$ at $x$ has the following geometric meaning. Suppose $x\in S^n$
is a regular point of $f$ and consider the restriction
$(\ker Df(x))^\perp\to\R^s$
of $Df(x)$ to the orthogonal complement
of $\ker Df(x)$.
Then $\NJ_f(x)$ of $f$ equals the absolute value of
the determinant of this map, cf.~(\ref{eq:volvol}).
The normal Jacobian of the differentiable map of Riemannian manifolds
$$
 F\colon S^n\times B_a \to \R^s, (x,t)\mapsto f(x) .
$$
at $(x,t)\in S^n\times B_a$ satisfies $\NJ_F(x,t) = \NJ_f(x)$.

Consider the following integrable function $\varphi_d$ for $\d>0$
$$
\varphi_\d\colon S^n\times B_a \to\R,\ \varphi(x,t)
 = \left\{\begin{array}{ll} 1 \ \mbox{ if $\|f(x)\|_\infty < \d$}\\
                            0 \ \mbox{ otherwise.}\end{array}\right.
$$
The coarea formula (cf.~\cite[Appendix]{howa:93} or \cite[III.\S2]{suwi:72})
applied to $F$ yields
\begin{gather*}
 \int_{y\in (-\d,\d)^s}  G_f(y)\ dy  
= \int_{y\in\R^s} \int_{f^{-1}(y)\times B_a} \varphi_\d\, g_f\, d(f^{-1}(y)\times B_a)\ dy \\
 = \int_{S^n\times B_a} \varphi_\d\, g_f\, \NJ_F\ d(S^n\times B_a) .
 \end{gather*}
Dividing by $(2\d)^s$ and taking the expectation over~$f$
with respect to the given Gaussian distribution, we obtain
\begin{gather}\label{eq:cabigeq} 
 \frac1{(2\d)^s} \int_{y\in (-\d,\d)^s} \E_f ( G_f(y) )\, dy 
 = \int_{S^n\times B_a} \frac{1}{(2\d)^s}\ \E_f \big(\varphi_\d\, g_f\, \NJ_F \big)\ d(S^n\times B_a).
\end{gather}
For fixed $(x,t)\in S^n\times B_a$ we can write the integrand $I_\d(x,t)$ 
on the right-hand side of (\ref{eq:cabigeq}) 
as an integral over conditional expectations as follows
\begin{equation}\label{eq:bdd}
 I_\d(x,t) = \frac{1}{(2\d)^s} \int_{y\in (-\d,\d)^s}
       p_{f(x)}(y)\, \E_f \left( g_f(x,t) \NJ_f (x)\Big/ f(x)=y \right)\, dy .
\end{equation}
By continuity, we get 
\begin{equation}\label{eq:limmy}
   \lim_{\d\to 0} I_\d(x,t) = p_{f(x)}(0)\,  \E_f \Big( g_f(x,t) \NJ_f(x) \Big/ f(x)=0\Big) .
\end{equation}
The integrand of (\ref{eq:bdd}) is a continuous function of~$(x,t,y)$
and therefore bounded by some constant $M$ on
$S^n\times B_a\times \{y\in\R^s\mid \|y\|\le 1\}$.
Hence (\ref{eq:bdd}) is as well bounded by~$M$ for all $0<\d\le 1$.
We may therefore apply Lebesgue's Theorem and interchange in (\ref{eq:cabigeq})
the integral over $(x,t)\in S^n\times B_a$ and the limit for $\d\to 0$ obtaining
\begin{multline}\label{eq:biggy}
\E_f(G_f(0)) = \lim_{\d\to 0} \frac1{(2\d)^s} \int_{y\in (-\d,\d)^s} \E_f ( G_f(y) )\, dy 
 =  \lim_{\d\to 0} \int_{S^n\times B_a} I_\d\, d(S^n\times B_a) \\
 = \int_{(x,t)\in S^n\times B_a} p_{f(x)}(0)\, \E_f \Big( g_f(x,t) \NJ_f(x)\Big/ f(x)=0 \Big)\ d(S^n\times B_a) ,
\end{multline} 
where we have used (\ref{eq:limmy}) for the last equality.

Note that for $r>0$ (with uniform random $u\in S^{s-1}$ independent of~$f$)
$$
 p_{f(x)}(0)\, \E_{f,u} \Big( g_f(x,ru) \NJ_f(x)\Big/ f(x)=0\Big)  
 = \frac1{(1+r^2)^{(n+1)/2}}\psi(r,x).
$$
Hence, using polar coordinates $t=ru$, the right-hand side of~(\ref{eq:biggy}) can be written as
$$
 \Oh_n\Oh_{s-1} \int_{0}^a \frac{r^{s-1}\Psi(r)}{(1+r^2)^{(n+1)/2}}\, dr .
$$
Taking into account (\ref{eq:globy}), this completes the proof.
\end{proof}

\subsection{Expected characteristic polynomial of Weingarten map}\label{se:exp-charpol}

In order to prove the main Theorem~\ref{th:main}
we will evaluate Theorem~\ref{th:rice} for independent Gaussian polynomials~$f_\s$ having
invariant distributions. We write $\d_\s =\d(f_\s)$ for the parameter of $f_\s$.
We may assume without loss of generality that $\E\,f_\s(q)^2 =1$
for all $1\le\s\le s$
(scaling does not change the parameter of $f_\s$).
Hence $f_\s(q)$ is standard normal and the joint distribution
of $f(q)$ has the density
$p_{f(q)}(y)=(2\pi)^{-s/2}\exp(-\frac12(y_1^2+\cdots y_s^2))$.
In particular, $p_{f(q)}(0)=(2\pi)^{-s/2}$.

We proceed by a sequence of intermediate steps.
For $(r,u)\in [0,\infty)\times S^{s-1}$ and
a fixed matrix $M\in\R^{s\times n}$ of rank~$s$
we consider the following conditional expectation
\begin{equation}
\calE (r,u,M):= \E_{f} \left(\det(\id - r Lf(q,u))\Big/ f(q)=0, Df(q)=M \right).
\end{equation}
Recall the characterization~(\ref{eq:defL}) of $Lf(q,u)$ in terms of $Df(q)$ and $D^2f(q)$.
From Lemma~\ref{le:dervq-invar} we know that $D^2f(q)$ is independent of
$Df(q)$. Hence the above expectation may be taken with respect to the distribution
of $D^2f(q)$ conditioned solely on the event $f(q)=0$.

\begin{lemma}\label{le:evalcalE}
For $(r,u)\in [0,\infty)\times S^{s-1}$ and
a fixed matrix $Df(q)$ of rank~$s$ we set
\begin{equation}\label{eq:def-V}
  v=(v_1,\ldots,v_s)^T := \diag(\d_1^{1/2},\ldots,\d_s^{1/2})\,  (Df(q)^\dagger)^T Qf(q)(u) .
\end{equation}
Then we have (recall (\ref{eq:def-gamma})
$$
 \calE (r,u,Df(q)) = \sum_{j=0}^{\lfloor\frac{n-s}{2}\rfloor} r^{2j}\gamma_{2j} {n-s\choose 2j}
   \bigg(\sum_{\s=1}^s (1-\delta_\s) v_\s^2 \bigg)^j .
$$
\end{lemma}

\begin{proof}
Proposition~\ref{pro:expdet} is the key to this result.
As usual let $T_qZ$ denote the kernel of $Df(q)$ and recall that
the Hessian~$Hf(q)$ was defined as
the restriction of the bilinear map $D^2f(q)$ to $T_qZ\times T_qZ$.
For the following introduce an orthonormal  basis adapted to $\R^n=T_qZ\oplus (T_qT)^\perp$
(or observe that~\S\ref{se:invar-mat}
could have beeen presented in a coordinate-free way).

The matrix of $Lf(q,u)$ is $O(n-s)$-invariant and Gaussian.
The same is true for the random matrix
$Lf(q,u)_{\mathrm{cond}}$, which is defined as the random matrix $Lf(q,u)$ conditioned on $f(q)=0$.
In order to apply Proposition~\ref{pro:expdet} we need to calculate
the parameter of $Lf(q,u)_{\mathrm{cond}}$.

We write $Hf(q)=(Hf_1(q),\ldots,Hf_s(q))$.
By identifying (bi)linear maps with their matrices we obtain from
Equation (\ref{eq:defL}) that
$$
 Lf(q,u) = -\sum_{\s=1}^s \frac{v_\s}{\sqrt{\d_\s}}\, Hf_\s(q) ,
$$
where
$v=(v_1,\ldots,v_s)$ is defined as in (\ref{eq:def-V}).
Hence Lemma~\ref{le:suppositionW} implies
$$
\delta_L:=\delta(Lf(q,u)_{\mathrm{cond}})= \sum_{\s=1}^s \frac{v_\s^2}{\d_\s} \, \delta(Hf_\s(q)_{\mathrm{cond}})
$$
using obvious notation.
Lemma~\ref{le:restrict-invar-distr} tells us that
the restriction of $f_\s$ to $T_qZ$ (whose distribution is invariant under the
orthogonal group of $T_qZ$) has the same parameter~$\delta_\s$ as $f_\s$.
Corollary~\ref{cor:delta-cond-appl} gives that
$\delta(Hf_\s(q)_{\mathrm{cond}})=\delta_\s(1-\delta_\s)$, hence
$\delta_L = \sum_\s (1-\delta_\s)v_\s^2$.
Proposition~\ref{pro:expdet} implies now with $X\sim N(0,1)$
$$
 \calE (r,u,Df(q)) = \E\, (1 + r\,\sqrt{\delta_L}\, X)^{n-s} .
$$
Hence, taking into account that $\gamma_{2j} = \E\, X^{2j}$, we conclude
$$
  \E\, (1 + r\,\sqrt{\delta_L}\, X)^{n-s} =
  \sum_{j=0}^{\lfloor\frac{n-s}{2}\rfloor}{n-s\choose 2j}\, r^{2j}\,\delta_L^j\,\gamma_{2j},
$$
which proves the lemma.
\end{proof}

\subsection{Proof of main theorem}\label{se:proofofmain}

We prove now the following reformulation of Theorem~\ref{th:main}
for zero sets in spheres.

\begin{theorem}\label{th:hauptsatz1}
Suppose that $f_\s\in H_{d_\s,n}$ are independent centered random polynomials with
$O(n+1)$-invariant Gaussian distribution of parameter $\delta_\s$.
Consider the random zero set $\mZ(f)\subseteq S^n$ where $f=(f_1,\ldots,f_s)$.
Then the expectation of the curvature coefficient $K_{s+2j}(\mZ(f))$ satisfies
$$
\frac{\E\,K_{s+2j}(\mZ(f))}{\Oh_{n-s-2j}\Oh_{s+2j-1} } = (\delta_1\cdots\delta_s)^{1/2}\,
  \sum_{\nu\in\N^s,\,|\nu|=j}  (1-\delta_1)^{\nu_1}\cdots (1-\delta_s)^{\nu_s}
  C^{(1)}_{\nu_1} \cdots C^{(1)}_{\nu_s}
 $$
 for $0\le j\le\lfloor \frac{n-s}{2}\rfloor$, where
the coefficients $C^{(1)}_{k}$ can be characterized by the generating function
\begin{equation}\label{eq:C1genf}
 (1-Y)^{-1/2} = \sum_{k=0}^\infty C_k^{(1)}\, Y^k
  = 1 + \frac12 Y+ \frac38 Y^2+\frac5{16}Y^3 +\cdots .
\end{equation}
We have $C^{(1)}_0=1$ and
$C^{(1)}_{k} = \frac{1\cdot 3\cdot 5\cdots (2k-1)}{k!\, 2^k} = \frac{(2k)!}{4^k\,k!^2}$
for $k>0$.

In the case where all $\delta_\s:=\delta$ are equal the result simplifies to
$$
 \E\,K_{s+2j}(\mZ(f)) = \delta^{s/2}\,(1-\delta)^j\, \Oh_{n-s-2j}\Oh_{s+2j-1} C_j^{(s)}\quad
  \mbox{ for $0\le j\le\lfloor \frac{n-s}{2}\rfloor$,}
$$
where the $C_j^{(s)}$ are characterized as
the coefficients of the power series
$$
(1-Y)^{-s/2} =\sum_{k=0}^\infty C_k^{(s)}\, Y^k .
$$
More specifically, we have $C_0^{(s)}=1$ and for $k>0$
$$
C_k^{(s)} = \frac{s(s+2)(s+4)\cdots (s+2k-2)}{k!\,2^k} .
$$
\end{theorem}

\proofof{Theorem~\ref{th:hauptsatz1}}
Put $\Delta:=\diag(\sqrt{\d_1},\ldots,\sqrt{\d_s})$.
Since $\d_\s^{-1/2}Df_\s(q)$ is standard normal distributed in $\R^n$ and
the $f_\s$ are independent, we can write $Df(q) = \Delta A$,
where $A\in\R^{s\times n}$ is a random matrix with independent
standard Gaussian entries. Note that
$\det(Df(q)Df(q)^T) = \d_1\cdots \d_s\, \det (A A^T)$
and
$Df(q)^\dagger = A^\dagger \Delta^{-1}$.

Let $u\in S^{s-1}$ be a uniformly distributed
random unit vector which is independent of $A$.
The function $\psi$ introduced in Theorem~\ref{th:rice} satisfies for $r>0$
$$
 \psi(r,q) = (2\pi)^{-s/2} (\d_1\cdots \d_s)^{1/2}\,
 \E_{A,u} \left(\sqrt{\det A A^T}\ \calE(r,u, \Delta A) \right)  .
$$
Lemma~\ref{le:evalcalE} tells us that
$$
 \calE (r,u, \Delta A)) = \sum_{j=0}^{\lfloor\frac{n-s}{2}\rfloor} r^{2j}\gamma_{2j} {n-s\choose 2j}
   \bigg(\sum_{\s=1}^s (1-\delta_\s) v_\s^2 \bigg)^j ,
$$
where $v\in\R^s$ is the image of $u$ under the linear endomorphism
$$
 \Delta (Df(q)^\dagger)^T Qf(q)(u)  = (A^\dagger)^T Q_A(u)
$$
of $\R^s$ (recall the definition of $Q_A$ in~(\ref{eq:QA})).
We make the multinomial expansion
$$
   \bigg(\sum_{\s=1}^s (1-\delta_\s) v_s^2 \bigg)^j
   = \sum_{\nu\in\N^s,\,|\nu|=j}{j\choose\nu}(1-\delta_1)^{\nu_1}\cdots (1-\delta_s)^{\nu_s}
         v_1^{2\nu_1}\cdots v_s^{2\nu_s} .
$$
Thus we need to compute for $\nu\in\N^s$ with $|\nu|=j$
$$
   \E_{A,u} \left(\sqrt{\det A A^T}\  v_1^{2\nu_1}\cdots v_s^{2\nu_s} \right)  .
$$
Proposition~\ref{le:heureka} determines the joint distribution of
$(\sqrt{\det A A^T},  v)$
for random~$A$ and $u$.
Accordingly, we write with a uniformly distributed $w\in S^{s-1}$ that is
independent of~$A$:
$$
   \E_{A,u} \left(\sqrt{\det A A^T}\  v_1^{2\nu_1}\cdots v_s^{2\nu_s} \right)
   = \E_{A,u}\left( \vol(A_1,\ldots,A_{s-1})\,
      \frac1{\|A_s^\perp\|^{2j-1}} w_1^{2\nu_1}\cdots w_s^{2\nu_s}  \right).
$$
It is well known that~\cite{weyl:39}
$$
 \E_{w\in S^{s-1}}\left( w_1^{2\nu_1}\cdots w_s^{2\nu_s}  \right) =
  \frac{\gamma_{2\nu_1}\cdots\gamma_{2\nu_s} }{s(s+2)\cdots (s+2j-2)} .
$$
By using Lemma~\ref{le:techcalc} we obtain
$$
  \E_{A,u} \left(\sqrt{\det A A^T}\  v_1^{2\nu_1}\cdots v_s^{2\nu_s} \right)
 = \frac{\Oh_{n-s}}{2 (2\pi)^{\frac{n-s}{2}}} \cdot
 \frac{\gamma_n\gamma_{n-s+1-2j}}{\gamma_{n-s+1}}\,
 \frac{\gamma_{2\nu_1}\cdots\gamma_{2\nu_s}}{s(s+2)\cdots (s+2j-2)} .
$$
This formula can be considerably simplified. We put
$C^{(1)}_{k}:=\frac{(2k)!}{4^k\, k!^2}$ for $k\in\N$.

\medskip

\noindent{\bf Claim.} We have
\begin{equation*}\label{eq:Ohgamma}
 \frac{\Oh_n\Oh_{s-1}}{(2\pi)^{s/2}\Oh_{n-s-2j}\Oh_{s+2j-1}}
   \gamma_{2j} {n-s\choose 2j} {j \choose \nu}
  \E_{A} \left(\sqrt{\det A A^T}\  v_1^{2\nu_1}\cdots v_s^{2\nu_s} \right)
  = C^{(1)}_{\nu_1}\cdots C^{(1)}_{\nu_s}.
\end{equation*}
\smallskip
In order to verify this recall first that
$$
\Oh_n = \frac{2\pi^{\frac{n+1}{2}}}{\Gamma(\frac{n+1}{2})},\quad
 \gamma_n = \frac1{\sqrt{\pi}}\, 2^{n/2}\,\Gamma(\frac{n+1}{2}) .
$$
From $\Gamma(x+1)=x\Gamma(x)$ and $\Gamma(1/2)=\sqrt{\pi}$ we get
$$
 \Gamma(m+1)=m!,\quad \Gamma(m+1/2)= (m-1/2)(m-3/2)\cdots 1/2 \sqrt{\pi}.
$$
Using the above, it is straightforward to check that
\begin{equation*}\label{eq:uno}
 \Oh_n\,\gamma_n = 2 (2\pi)^{n/2},\quad
 \frac{n!\,\Oh_n}{\gamma_{n+1}} = (2\pi)^{\frac{n+1}{2}},\quad
  (2\pi)^j\,\frac{\Oh_{s-1}}{\Oh_{s-1+2j}} = s(s+2)\cdots (s+2j-2).
\end{equation*}
Moreover, recall from~(\ref{eq:gamma}) that
$(2j)! = 2^j\, j! \gamma_{2j}$ and note that
$C_k^{(1)}=\frac{\gamma_{2k}}{2^k\,k!}$.
The claim follows by simplifying the formula using the above stated equations
in a straightforward (but tedious) way.

\medskip

Combining what we have shown so far we obtain
$$
 \frac{\Oh_n\Oh_{s-1}\, \psi(r,q)}{\Oh_{n-s-2j}\Oh_{s+2j-1}} =
  (\d_1\cdots \d_s)^{1/2} \sum_{|\nu| \le \lfloor \frac{n-s}{2}\rfloor}
   (1-\delta_1)^{\nu_1}\cdots (1-\delta_s)^{\nu_s}
      C^{(1)}_{\nu_1}\cdots C^{(1)}_{\nu_s}  r^{2|\nu|}.
$$

We apply now Theorem~\ref{th:rice}.
From the $O(n+1)$-invariance it follows that
$\psi(r,x)=\psi(r,q)$ for all $x\in S^n$,
hence $\Psi(r) = \psi(r,q)$.
Let $0<\a<\pi/2$ and put $a=\tan\a$.
By substituting $r=\tan\rho$ we obtain
\begin{equation*}\label{eq:Jdef}
 \int_0^a \frac{r^{s+2j-1}}{(1+r^2)^{(n+1)/2}}\,dr = J_{n,s+2j}(\a) .
\end{equation*}
We conclude with Theorem~\ref{th:rice} that
\begin{eqnarray*}
 \lefteqn{\sum_{j=0}^{\lfloor\frac{n-s}{2}\rfloor}
       \frac{\E\,K_{s+2j}(\mZ(f))}{\Oh_{n-s-2j}\Oh_{s+2j-1}}\, J_{n,s+2j}(\a)}\\
  &=& \frac{\Oh_n\Oh_{s-1}}{\Oh_{n-s-2j}\Oh_{s+2j-1}}
       \int_{0}^a \frac{r^{s-1}\,\Psi(r)}{(1+r^2)^{(n+1)/2}}\, dr \\[2ex]
  &=& (\d_1\cdots \d_s)^{1/2} \sum_{|\nu| \le \lfloor \frac{n-s}{2}\rfloor}
   (1-\delta_1)^{\nu_1}\cdots (1-\delta_s)^{\nu_s}
      C^{(1)}_{\nu_1}\cdots C^{(1)}_{\nu_s}\, J_{n,s+2|\nu|}(\a) .
\end{eqnarray*}
By comparing the coefficients of the linearly independent functions $J_{n,k}(\a)$,
the stated formula for $\E\, K_{s+2j}(\mZ(f))$ follows.

To settle the case where all $\d_\s$ are equal just note
$(1-Y)^{-1/2} =\sum_{k=0}^\infty C_k^{(1)}\, Y^k$ and
$(1-Y)^{-s/2} =\sum_{k=0}^\infty C_k^{(s)}\, Y^k$
implies that
$C^{(s)}_j = \sum_{|\nu|=j} C^{(1)}_{\nu_1}\cdots C^{(1)}_{\nu_s}$.
\proofend

\section{Alternative proof of the main result}\label{se:comm}

We show here that Theorem~\ref{th:main} can be quickly derived from the knowledge
of the expected Euler characteristic of a random projective hypersurface $\mZ(f)$
for an invariant centered Gaussian random polynomial~$f$.
The key of this reduction is the kinematic formula and the generalized Gauss-Bonnet theorem.

In this section, $\mZ(f)$ stands for the zero set in $\proj^n$.
Suppose $f\in H_{d,2\ell+1}$ has invariant centered Gaussian distribution with parameter $\d$.
We define $\chi_\ell(\d) := \E\,\chi(\mZ(f))$.
Theorem~\ref{th:main} in the case of one equation ($s=1$) yields
\begin{equation}\label{eq:chi-one-equ}
 \chi_\ell(\d) = \d^{1/2}\sum_{k=0}^{\ell} C_k^{(1)} (1-\d)^{k} .
\end{equation}
Taking into account Equation (\ref{eq:C1genf}),
we get a closed form expression for the generating function $\chi(\d;T)$
of $\chi_\ell(\d)$ as follows:
\begin{eqnarray}\notag
\chi(\d;T) &:=& \sum_{\ell=0}^\infty \chi_\ell(\d)\, T^{2\ell}
   = \d^{1/2}\sum_{k=0}^\infty C_k^{(1)} (1-\d)^{k}  T^{2k} \sum_{\ell=k}^{\infty}  T^{2(\ell-k)} \\ \label{eq:chiseries}
  &=& \frac{\d^{1/2}}{(1-T^2)(1-(1-\d)T^2)^{1/2}}.
\end{eqnarray}
This formula can also be readily deduced from
Podkorytov's result~\cite{podk:99}, cf.~(\ref{eq:podk}).

Using the kinematic formula we can prove a stability result for
$\E\, \mu_e(\mZ(f))$.

\begin{lemma}\label{le:restrict-invar-distr1}
Suppose $f\in H_{d,n}$ is $O(n+1)$-invariant and $n'\le n$.
Then the restriction $f'\in H_{d,n'}$ of $f$ to $\R^{n'+1}$ is $O(n'+1)$-invariant
and has the same parameter.
For $0\le e < n'$, $e$ even, we have
$\E\,\mu_e(\mZ(f')) = \E\, \mu_e(\mZ(f))$.
\end{lemma}

\begin{proof}
The first part of the statement was already established in Lemma~\ref{le:restrict-invar-distr}.
By Theorem~\ref{th:kinform-S} (in the version for $\proj^n$)
we have for almost all $f\in H_{d,n}$, as $\mu(\proj^{n'};T)=1$,
$$
 \int \mu(\mZ(f)\cap g\proj^{n'};T)\, dg
 \equiv \mu(\mZ(f);T) \bmod T^{n'},
$$
where the integral is with respect to the Haar measure of $O(n+1)$ scaled such that
the volume of $O(n+1)$ equals~$1$.
Taking the expectation over $f$ and interchanging with the integral over $g$
we obtain
$$
 \int
 \E_f(\mu(\mZ(f)\cap g\proj^{n'};T))\, dg
 \equiv \E_f\,\mu(\mZ(f);T)\bmod T^{n'}.
$$
By the invariance of the distribution of $f$, the integrand is independent of $g$,
hence the integral equals $\E_f(\mu(\mZ(f)\cap \proj^{n'};T))$.
This expectation equals
$\E_{f'}(\mu(\mZ(f');T))$, which finishes the proof.
\end{proof}

\begin{lemma}\label{le:muchi}
Suppose that $f\in H_{d,n}$ is $O(n+1)$-invariant with parameter~$\d$.
Then we have $\E\,\mu_{0}(\mZ(f)) = \chi_0(r)$
and for $1\le k\le(n-1)/2$
$$
 \E\,\mu_{2k}(\mZ(f)) = \chi_k(\d) - \chi_{k-1}(\d).
$$
In particular, $\E\,\mu_{2k}(\mZ(f))$ depends only on $k$ and $\d$
and not on the dimension~$n$ of the ambient space.
\end{lemma}

\begin{proof}
We denote by $f'$ and $f''$ the restrictions of $f$ to $\R^{2k+2}$
and $\R^{2k}$, respectively.
Theorem~\ref{th:GB-S} (in the version for $\proj^n$) implies that
$$
 \chi_k(\d)     = \sum_{i=0}^k \E\, \mu_{2i}(\mZ(f')),\quad
 \chi_{k-1}(\d) = \sum_{i=0}^{k-1} \E\, \mu_{2i}(\mZ(f'')).
$$
By Lemma~\ref{le:restrict-invar-distr1} we have
$\E\,\mu_{2i}(\mZ(f)) = \E\,\mu_{2i}(\mZ(f')) = \E\,\mu_{2i}(\mZ(f''))$
for $0\le i\le k-1$ and
$\E\,\mu_{2k}(\mZ(f)) = \E\,\mu_{2k}(\mZ(f'))$.
Subtracting the above two equations, we obtain
$\chi_k(\d)-\chi_{k-1}(\d) = \E\,\mu_{2k}(\mZ(f))$
as claimed. The assertion for $k=0$ is obvious.
\end{proof}

We proceed now with an alternative proof of Theorem~\ref{th:main}

\proofof{Theorem~\ref{th:main}}
By Lemma~\ref{le:muchi} the formal power series
$$
 \mu(\d;T):= \sum_{k=0}^\infty \E\,\mu_{2k}(\mZ(f))\, T^{2k}
$$
satisfies $\mu(\d;T)= (1-T^2)\,\chi(\d;T)$.
Equation~(\ref{eq:chiseries}) implies that
\begin{equation}\label{eq:muv}
 \E\,\mu(\d;T) = \d^{1/2}(1-(1-\d)T^2)^{-1/2}.
\end{equation}
Suppose now that $f_\s\in H_{d_\s,n}$ are independent random variables
with invariant centered Gaussian distribution of parameter $\d_\s$
for  $1\le \s\le s\le n$. Write $\underline{f} :=(f_1,\ldots,f_{s-1})$.
The kinematic formula (Theorem~\ref{th:kinform-S} in the version for $\proj^n$)
tells us that for almost all $\underline{f},f_s$
$$
 \int \mu(\mZ(\underline{f})\cap g\mZ(f_s);T)\, dg
 \equiv \mu(\mZ(\underline{f});T)\ \mu(\mZ(f_s);T)\bmod T^{n-s+1},
$$
where the integral is over $O(n+1)$ with respect to the Haar measure scaled to~$1$.
We take the expectation with respect to $\underline{f}$ and $f_s$.
Taking their independence into account we get
$$
 \int \E_{\underline{f},f_s}\,\mu(\mZ(\underline{f})\cap g\mZ(f_s);T)\, dg
 \equiv \E\,\mu(\mZ(\underline{f});T)\ \E\,\mu(\mZ(f_s);T) \bmod T^{n-s+1}.
$$
By invariance, the integrand does not depend on $g$ and equals
$\E\,\mu(\mZ(f_1,\ldots,f_s);T)$.
We conclude by induction that
$$
\E\,\mu(\mZ(f_1,\ldots,f_s);T)\equiv
 \E\,\mu(\mZ(f_1);T)\cdots \E\,\mu(\mZ(f_s);T) \bmod T^{n-s+1}.
$$
Plugging in the explicit expression (\ref{eq:muv}) for $\E\,\mu(\mZ(f_i);T)$,
the desired formula for the expectation of the curvature polynomial follows.
Finally note that
$ (1-Y)^{-s/2} = \sum_{k=0}^\infty C_k^{(s)} Y^k $.
\endproof

\begin{remark}\label{re:nogauss-assumpt}
The above reduction to the computation of the expected
Euler characteristic works for {\em any}  $O(n+1)$-invariant distribution of random polynomials.
The Gaussian assumption is not needed for the reduction.
\end{remark}

\begin{remark}\label{re:adta}
We briefly outline how (\ref{eq:chiseries}) can be
derived from a general result of Taylor and Adler~\cite{taad:03}.
Suppose that $f\in H_{d,n}$ is $O(n+1)$-invariant with parameter~$\d$
and suppose w.l.o.g.\ that $\E\,f(p)^2 =1$ for $p\in S^n$.
We consider $f$ as a centered  unit variance Gaussian field on the sphere $S^n$.
Then $f$ defines a Riemannian metric on~$S^n$ by
$g_p(X_p,Y_p) := \E ( X_p f\cdot Y_p f)$ for tangent vectors $X_p,Y_p\in T_p S^n$.
It easily follows that $g_p(X_p,Y_p) = \d\, \langle X_p,Y_p \rangle$,
where $\langle\ ,\ \rangle$ denotes the scalar product on $T_pS^n$.
Hence, with respect to this metric, $S^n$ is isometric to the sphere $M$ in
$\R^{n+1}$ of radius $R=\sqrt{\d}$ .
Theorem 4.1 of \cite{taad:03} states that
$$
 \E\, \chi( S^n \cap f^{-1}[u,\infty)) = \sum_{j=0}^n \mathcal{L}_j \rho_j(u) ,
$$
with ``Lipschitz-Killing curvatures'' $\mathcal{L}_{j}$ of $M$ and
functions $\rho_j(u) $ related to Hermite polynomials.
Almost surely, $S^n \cap f^{-1}[u,\infty)$ is a compact domain with smooth
boundary $S^n\cap f^{-1}(0)$, which implies
$\chi(S^n\cap f^{-1}(0)) = 2 \chi( S^n \cap f^{-1}[u,\infty))$ if $n$ is odd.
The Lipschitz-Killing curvatures of $M$ can be defined via Weyl's formula for
the volume of the tubes around $M$ in $\R^{n+1}$
$$
 \vol(T(M,r)) =  \sum_{j=0}^n \mathcal{L}_{j}\, \omega_{n+1-j}\, r^{n+1-j}
$$
with $\omega_0:=1$ and $\omega_k= \Oh_{k-1}/k$ for $k>0$.
A straightforward calculation yields
$\mathcal{L}_j =\frac{2\omega_{n+1}}{\omega_{n+1-j}}\, {n+1\choose j} R^j$
if $j\equiv n\bmod 2$ and $\mathcal{L}_j =0$ otherwise.
Hence we obtain for odd~$n$
$$
 \E\,\chi(\mZ_{\proj^n}(f)) = \frac12 \E\, \chi( S^n \cap f^{-1}(0))
  = \sum_{1\le j\le n, \mbox{\scriptsize{$j$ odd}}} \frac{2\omega_{n+1}}{\omega_{n+1-j}}\,
     {n+1\choose j} \d^{j/2} \rho_j(0) .
$$
It is possible to derive Equation~(\ref{eq:chi-one-equ}) from this, but we omit the details.
\end{remark}

\section*{Appendix}\label{se:appendix}

\proofof{Theorem~\ref{th:GB-S}}
W.l.o.g.\ $m<n$.
We suppose first that $n$ is odd.
The submanifold $M$ is a deformation retract of the tube $T_\a:=T(M,\a)$
for sufficiently small $\a>0$,
hence $\chi(M)=\lim_{\a\to 0}\chi(T_\a)$.
The generalized Gauss-Bonnet formula
applied to the domain $T_\a$ in $S^n$ says that
(cf.~\cite[(17.22), p.~303]{sant:76})
$$
 \frac12\Oh_n\,\chi(T_\a) =
 \sum_{0\le i \le n-1,\,\mbox{\scriptsize $i$ even}} \frac{\Oh_n}{\Oh_{n-1-i}\Oh_i}\
  \int_{\partial T_\a} \sigma_i(\kappa_1,\ldots,\kappa_{n-1}) d (\partial T_\a) .
$$
Going to the limit $\a\to 0$ and applying
Equation~(\ref{eq:meancurv}) we get
($s=n-m$ is odd)
\begin{eqnarray*}
 \frac12\chi(M) &=&
  \sum_{0\le i \le n-1,\,\mbox{\scriptsize $i$ even}}
  \frac{1}{\Oh_{n-1-i}\Oh_i}\, K_{i+1}(M) \quad \mbox{\small (put $i+1=s+e$)}\\
  &=&\sum_{0\le e \le m,\,\mbox{\scriptsize $e$ even}}
  \frac{1}{\Oh_{m-e}\Oh_{s+e-1}}\, K_{s+e}(M)
   = \sum_{0\le e \le m,\,\mbox{\scriptsize $e$ even}} \mu_e(M).
\end{eqnarray*}

In the case where $n$ is even, we argue similarly:
The generalized Gauss-Bonnet formula (cf.~\cite[(17.21), p.~303]{sant:76}
applied to $T_\a$ says that
$$
 \frac12\Oh_n\chi(T_\a) = \vol(T_\a) +
 \sum_{1\le i \le n-1,\,\mbox{\scriptsize $i$ odd}} \frac{\Oh_n}{\Oh_{n-1-i}\Oh_i}\
  \int_{\partial T_\a} \sigma_i(\kappa_1,\ldots,\kappa_{n-1}) d(\partial T_\a).
$$
Hence, by taking the limit $\a\to 0$
and using~(\ref{eq:meancurv}), we get (note that $s$ is even)
\begin{eqnarray*}
 \frac12\chi(M) &=&
  \sum_{1\le i \le n-1,\,\mbox{\scriptsize $i$ odd}}
  \frac{1}{\Oh_{n-1-i}\Oh_i} \,K_{i+1}(M) \\
   &=& \sum_{0\le e \le m,\,\mbox{\scriptsize $e$ even}}
       \frac{1}{\Oh_{m-e}\Oh_{s+e-1}} \,K_{s+e}(M)
          = \sum_{0\le e \le m,\,\mbox{\scriptsize $e$ even}} \mu_e(M),
\end{eqnarray*}
which finishes the proof.
\proofend

\proofof{Proposition~\ref{pro:delta=}}
The covariance function
$r(x,y):=\E(f(x)f(y))$
is a polynomial function that is homogeneous of degree~$d$ in both
sets of variables $x\in\R^{n+1}$ and $y\in\R^{n+1}$.
Since the distribution of $f$ is invariant under the action of the orthogonal group,
we have $r(gx,gy)=r(x,y)$ for all $g\in O(n+1)$ and $x,y\in\R^{n+1}$.
It follows from invariant theory that $r(x,y)$ is a real polynomial in $\|x\|^2,\|y\|^2,\langle x,y\rangle$
(e.g., see~\cite{spiv5:79}). From the fact that $r$ is bihomogeneous of degree $(d,d)$
it is easy to conclude that $r$ has the following form
$$
 r(x,y) = \sum_{k=0}^{\lfloor d/2\rfloor} \beta_k\, \|x\|^{2k} \|y\|^{2k}\langle x,y\rangle^{d-2k}
 \quad\quad (\beta_k\in\R) .
$$

We can express the moments of the partial derivatives of $f$ at $q$
by partial derivatives of the covariance function~$r$ at $(q,q)$:
we have for $1\le i,j,k,\ell\le n$
\begin{gather} \notag
 \E\, (\partial^2_{x_i x_j} f(q)\, f(q)) = \partial^2_{x_i x_j} r(q,q),\
 \E\, (\partial_{x_i} f(q)\, \partial_{y_k}  f(q)) = \partial^2_{x_i y_k} r(q,q), \\ \label{eq:rder}
 \E\, (\partial^2_{x_i x_j} f(q)\, \partial^2_{y_k y_\ell} f(q))
  = \partial^4_{x_i x_jy_k y_\ell} r(q,q) .
\end{gather}

In order to show the second claim it is convenient to use the abbreviation
$a_{ij}:= \partial^2_{x_i x_j}f(q)$, $A:=(a_{ij})_{1\le i,j\le n}$.
By Lemma~(\ref{le:2nd-deriv}) we have $W:=D^2f(q) = A -d f(q) I_n$.
We obtain for the parameter of $W$
\begin{eqnarray*}
 \lefteqn{\delta(W) = \frac1{n(n-1)}\Big( \E\,(\tr W)^2 - \E\,\|W\|_F^2 \Big)}\\
  & = & \frac1{n(n-1)}\Big( \E\,(\tr A)^2 - \E\,\|A\|_F^2
          - 2(n-1) d\, \E\,(f(q)\tr A) + n(n-1) d^2 \E\, f(q)^2 \Big)\\
  & = & \frac1{n(n-1)}\sum_{i\ne j} \E\,(a_{ii} a_{jj} - a_{ij}^2)
          -\frac{2}{n} d\, \E\,(f(q)\,\tr A) + d^2 \E\, f(q)^2 \\
  & = & \frac1{n(n-1)}\sum_{i\ne j}\Big(\partial^4_{x_i x_i y_j y_j}r
          - \partial^4_{x_i x_j y_i y_j}r\Big)(q,q)
          -\frac{2}{n} d\, \sum_i \partial^2_{x_i x_i}r(q,q) + d^2 r(q,q) .
\end{eqnarray*}
On the other hand,
$$
\delta(f)\E\,f(q)^2 = \frac1{n}\sum_i\E\;(\partial_{x_i}f(q))^2
  = \frac1{n} \sum_i \partial^2_{x_i y_i} r(q,q) .
$$
In order to prove the second claim, it suffices to check equality of
the above two expressions.
Since both expressions are linear in~$r$,
it is enough to check this for
$r_k(x,y)= \|x\|^{2k} \|y\|^{2k}\langle x,y\rangle^{d-2k}$.

A tedious but straightforward calculation  yields for
$1\le i,j\le n$, $i\ne j$,
\begin{gather*}
 r_k(q,q)=1,\
 \partial_{x_i x_i} r_k(q,q) = 2k,\
 \partial_{x_i y_i} r_k(q,q) = d-2k,\\
 \partial^4_{x_i x_i y_j y_j}r(q,q) = 4k^2,\
 \partial^4_{x_i x_j y_i y_j}r(q,q) = (d-2k)(d-2k-1).
\end{gather*}
Plugging in this in the above expressions we see that indeed
$\delta(W)= \delta(f)\E\, f(q)^2$.
The verification of the first claim is similar and a bit simpler.
\proofend

\proofof{Lemma~\ref{le:2ff-implicit}}
By invariance of the assertion under the orthogonal group it is sufficient to verify
the claim at the point $q :=(1,0\ldots,0)$ and we may also assume that
$T_qZ = \ker Df(q) = \R^{n-s}\times 0^s$.
Hence $\partial_{n-s+\t}f_\s(q)=0$ for $1\le\s,\t\le s$.

Our assumption $\mathrm{rank}\, Df(q) =s$ allows to apply the implicit
function theorem.
There is an open subset $U\subseteq\R^{n-s}$ containing the origin and
a differentiable map
$h\colon\R^{n-s}\supseteq U \to\R^s$ such that $h(0)=0$ and
$$
 \varphi\colon\R^{n-s}\supseteq U \to S^n,\
  u\mapsto \frac{(1,u,h(u))}{\sqrt{1+\|u\|^2 + \|h(u)\|^2}}
$$
is a local diffeomorphism of $U$ onto an open neighborhood of $q$ in $Z$.

We make the following useful convention on indices:
$\a,\b,\g$ run in the range $1,2,\ldots,n-s$ while
$\s,\t,\rho$ run in the range $1,2,\ldots,s$.
A straightforward calculation shows that 
$$
  \partial_\a \varphi_0(0) = 0,\quad
  \partial_\a \varphi_\g(0) = \delta_{\a\g}, \quad
  \partial_{\a} \varphi_{n-s+\s}(0) = \partial_{\a} h_\s (0) = 0,
$$
where the last equality follows from our
assumption $T_qZ = \R^{n-s}\times 0^s$.
A similar calculation yields 
$$
  \partial^2_{\a,\b} \varphi_{n-s+\s}(0) = \partial^2_{\a,\b} h_\s (0).
$$
By the definition~(\ref{eq:def-2ndff}) of the second fundamental form we obtain for
$V,W\in T_qZ=\R^{n-s}\times 0^s$ and $\nu\in (T_xZ)^\perp=0^{n-s}\times\R^s$
\begin{equation}\label{eq:zwff}
\zwff_Z(x)(V,W,\nu) = \sum_{\a,\b,\s} \nu_\s \partial^2_{\a,\b} h_{\s}(0)V_\a W_\b .
\end{equation}
It remains to express $\partial^2_{\a,\b}h_{\s}(0)$ by partial derivatives of $f$.

By taking the derivative of $f_\rho(1,u,h(u)) = 0$ with respect to $u_\a$ we obtain
$$
 \partial_{X_\a} f_\rho(1,u,h(u)) +
 \sum_\t \partial_{X_{n-s+\t}} f_\rho(1,u,h(u))\cdot \partial_\a h_\t(u) = 0 .
$$
Differentiating this with respect to $u_\b$ and taking into account that
$\partial_\a h_\t(0) = 0$ we obtain after a short calculation at $v=0$
$$
 \partial^2_{X_\a,X_\b} f_\rho(q) +
 \sum_\t \partial_{X_{n-s+\t}} f_\rho(q)\cdot \partial^2_{\a,\b} h_\t(0) = 0 .
$$
We may write this as
$
 \partial^2_{\a,\b} h_\s(0) = - \sum_\rho m_{\s,\rho} \partial^2_{X_\a,X_\b} f_\rho(q),
$
where $(m_{\s,\rho})$ denotes the matrix of $Df(q)^\dagger$.
Plugging this into (\ref{eq:zwff}) we get
\begin{eqnarray*}
\zwff_Z(x)(V,W,\nu) &=&
 -\sum_{\s} \nu_\s \sum_\rho m_{\s,\rho}
  \sum_{\a,\b} \partial^2_{X_\a,X_\b} f_{\rho}(q)V_\a W_\b \\
 &=& \langle \nu, Df(q)^\dagger Hf(q)(V,W)\rangle,
\end{eqnarray*}
which was to be shown.
\proofend

{\small

}
\end{document}